\documentclass[11pt, reqno]{article}

\usepackage{textcmds}  
\usepackage{amsmath, amssymb, amsfonts, amstext, verbatim, amsthm, mathrsfs}
\usepackage[mathcal]{eucal}
\usepackage{microtype}
\usepackage{fullpage}
\usepackage{stmaryrd}
\usepackage{graphicx}
\usepackage[all]{xy}
\usepackage[modulo]{lineno}
\usepackage[usenames]{color}
\usepackage{enumitem}
\usepackage{xspace}
\usepackage{authblk}

\usepackage{setspace}
\usepackage{aliascnt} 
\usepackage[colorlinks=true,linkcolor=blue,citecolor=blue,urlcolor=blue,citebordercolor={0 0 1},urlbordercolor={0 0 1},linkbordercolor={0 0 1}]{hyperref} 
\usepackage{cleveref}

\usepackage[shortalphabetic]{amsrefs} 

\theoremstyle{plain}

\newcommand{\refnewtheoremn}[4]{%
\newaliascnt{#1}{#2}
\newtheorem{#1}[#1]{#3}
\aliascntresetthe{#1}
\expandafter\providecommand\csname #1autorefname\endcsname{#4}}

\newcommand{\refnewtheorem}[3]{\refnewtheoremn{#1}{#2}{#3}{#3}}

\newtheorem{thm}{Theorem}[section]

\refnewtheorem{lem}{thm}{Lemma}
\refnewtheorem{claim}{thm}{Claim}
\refnewtheorem{cor}{thm}{Corollary}
\refnewtheorem{conj}{thm}{Conjecture}
\refnewtheorem{prop}{thm}{Proposition}
\refnewtheorem{quest}{thm}{Question}
\refnewtheorem{amplif}{thm}{Amplification}
\refnewtheorem{scholium}{thm}{Scholium}

\refnewtheorem{tablethm}{thm}{Table}
\refnewtheorem{assumption}{thm}{Assumption}
\refnewtheorem{conclusion}{thm}{Conclusion}
\refnewtheorem{problem}{thm}{Problem}

\theoremstyle{definition}
\refnewtheorem{rem}{thm}{Remark}
\refnewtheorem{defn}{thm}{Definition}
\refnewtheorem{notn}{thm}{Notation}
\refnewtheorem{ex}{thm}{Example}
\refnewtheorem{fact}{thm}{Fact}
\refnewtheorem{hypothesis}{thm}{Hypothesis}

{\end{minipage} \end{center}}

\newcommand{\bb}[1]{\mathbb{#1}}
\newcommand{\mc}[1]{\mathcal{#1}}
\newcommand{\mf}[1]{\mathfrak{#1}}

\renewcommand{\rm}[1]{\mathrm{#1}}

\usepackage{xparse}

\ExplSyntaxOn

\NewDocumentCommand{\definealphabet}{mmmm}
 {
  \int_step_inline:nnn { `#3 } { `#4 }
   {
    \cs_new_protected:cpx { #1 \char_generate:nn { ##1 }{ 11 } }
     {
      \exp_not:N #2 { \char_generate:nn { ##1 } { 11 } }
     }
   }
 }

\ExplSyntaxOff

\definealphabet{bb}{\bb}{A}{Z}
\definealphabet{mc}{\mc}{A}{Z}
\definealphabet{mf}{\mf}{A}{Z}
\definealphabet{mf}{\mf}{a}{z}

\DeclareMathOperator{\Fun}{Fun}
\newcommand{\exact}{{\mathrm ex}}
\DeclareMathOperator{\DCoh}{DCoh}

\DeclareMathOperator{\QC}{QC}

\DeclareMathOperator{\ICoh}{ICoh}
\DeclareMathOperator{\Perf}{Perf}
\DeclareMathOperator{\APerf}{APerf}
\DeclareMathOperator{\cofib}{cofib}
\DeclareMathOperator{\colim}{colim}
\DeclareMathOperator{\gr}{gr}
\DeclareMathOperator{\ev}{ev}
\DeclareMathOperator{\Sym}{Sym}
\DeclareMathOperator{\uHom}{\underline{Hom}}
\DeclareMathOperator{\Hom}{Hom}
\DeclareMathOperator{\rank}{rank}
\DeclareMathOperator{\Spec}{Spec}

\DeclareMathOperator{\NS}{NS}
\DeclareMathOperator{\tot}{tot}
\DeclareMathOperator{\fib}{fib}

\DeclareMathOperator{\Coh}{\underline{Coh}}
\DeclareMathOperator{\spfilt}{sf}
\DeclareMathOperator{\Filt}{Filt}
\DeclareMathOperator{\Map}{Map}
\DeclareMathOperator{\Grad}{Grad}
\DeclareMathOperator{\Cycles}{\mathfrak{C}_\ast}

\newcommand{\loc}{\mathrm{loc}}
\newcommand{\nil}{\mathrm{nil}}
\newcommand{\radj}[1]{{\mathrm \beta}^{\geq #1}}
\newcommand{\ladj}[1]{{\mathrm \beta}^{< #1}}
\DeclareMathOperator{\hwt}{wt_{\rm{max}}}
\DeclareMathOperator{\wt}{wt}
\newcommand{\Mod}{\operatorname{-Mod}}
\DeclareMathOperator{\Aut}{Aut}

\title{A categorical perspective on non-abelian localization}
\author{Daniel Halpern-Leistner}
\affil{Cornell University}
\date{September 2025}

\begin{document}

\maketitle

\begin{abstract}
In equivariant geometry, a localization (a.k.a., concentration) theorem is typically interpreted as a relationship between the equivariant geometry of a space with a group action and the geometry of its fixed locus. We take a different perspective, that of non-abelian localization: a localization theorem relates the geometry of an algebraic stack that is equipped with a $\Theta$-stratification to the geometry of the centers of this stratification. We establish a ``virtual'' $K$-theoretic non-abelian localization formula, meaning it applies to algebraic derived stacks with perfect cotangent complexes. We also establish a categorical upgrade of this theorem, by introducing a category of ``highest weight $K$-homology cycles'' with respect to the stratification, and relating the category of highest weight cycles on the stack to those on the centers of its $\Theta$-stratification. We apply these results to prove a universal wall-crossing formula, and establish a new finiteness theorem for the cohomology of tautological complexes on the stack of one-dimensional sheaves on an algebraic surface.
\end{abstract}

\tableofcontents

\section{Introduction}

Let $X$ be a smooth projective algebraic variety with a linearizable $\bbC^\ast$-action. The Atiyah-Bott localization theorem \cite{Atiyah1984TheCohomology} compares the cohomology of the fixed locus $X^{\bbC^\ast}$, whose connected components we denote $Z_1,Z_2,\ldots$, with the equivariant cohomology of $X$. There are at least three versions of the localization theorem, which we state in topological $K$-theory rather than in cohomology:\footnote{In fact, the $K$-theoretic version preceded the cohomological version. In the context of algebraic $K$-theory, results like \ref{I:localization_modules} are sometimes called concentration theorems to disambiguate them from Quillen's localization theorem, and \ref{I:localization_homology} is sometimes called the trace formula.}
\begin{enumerate}[label=(\Alph*)]
\item The restriction map $K^{i}_{\bbC^\ast}(X) \to K^{i}_{\bbC^\ast}(X^{\bbC^\ast})$ becomes an isomorphism after inverting an element of the ground ring $K^0_{\bbC^\ast}(\mathrm{pt})$;\label{I:localization_modules}
\item The identity in $K^0_{\bbC^\ast}(X)$ decomposes as $1_X = \sum_\alpha (\sigma_\alpha)_\ast \left( \frac{1_{Z_\alpha}}{e(\bbN_{Z_\alpha/X})} \right)$, where $\sigma_\alpha : Z_\alpha \hookrightarrow X$ is the inclusion, $\bbN_{Z_\alpha/X}$ is the normal bundle of $Z_\alpha$, and $e(-)$ denotes the Euler class; and\label{I:localization_cohomology}
\item The fundamental class in equivariant $K$-homology concentrates on the fixed loci $Z_\alpha$, giving a formula for the $K$-theoretic index of any $E \in K^0_{\bbC^\ast}(X)$,\label{I:localization_homology}
\begin{equation}\label{E:first_localization_formula}
    \chi(X,E) = \sum_{\alpha} \chi\left(Z_\alpha, \frac{E|_{Z_\alpha}}{e(\bbN_{Z_\alpha/X})}\right).
\end{equation}
\end{enumerate}
In order to interpret \ref{I:localization_cohomology} and \ref{I:localization_homology}, one must invert $1-q^n \in K^0_{\bbC^\ast}(\rm{pt}) \cong \bbZ[q^{\pm 1}]$, where $n$ is a common multiple of all weights appearing in $\bbN_{Z_\alpha/X}$ for some $\alpha$, so that all of the $e(\bbN_{Z_\alpha/X})$ become units. One interprets \eqref{E:first_localization_formula} as an identity in the localization $\bbZ[q^{\pm 1}]_{1-q^n}$, where the right-hand-side a posteriori lies in $\bbZ[q^{\pm 1}]$.

Results of this kind were originally proved in equivariant topological $K$-theory by Atiyah and Segal \cite{Atiyah1968TheII} for the action of an arbitrary compact group $G_{\rm{cpct}}$ on a manifold, and \ref{I:localization_modules} is proved in \cite{MR234452} for arbitrary locally compact $G_{\rm{cpct}}$-spaces. Since then, they have been extended in several directions. A version of \ref{I:localization_modules} holds in algebraic $K$-theory for actions of a split reductive group scheme over a noetherian base on a separated finite type algebraic space \cite[Thm.~2.2]{Thomason1992UneAlgebrique}, and a relative version of \ref{I:localization_homology} holds for the pushforward along a proper $T$-equivariant morphism $X \to Y$, where $T$ is a diagonalizable group scheme (e.g., a torus) and $X$ and $Y$ are finite type separated algebraic spaces \cite[Thm.~3.5]{Thomason1992UneAlgebrique}. The latter, however, is not canonical, as it makes use of a closed equivariant embedding $X \hookrightarrow Z$, where $Z$ is regular and proper over $Y$.

When $X$ is singular but has a perfect obstruction theory, which in practice occurs because $X$ is the classical scheme underlying a quasi-smooth derived scheme, a version of \ref{I:localization_homology} was established for the virtual fundamental class in \cite{Graber1999LocalizationClasses}. The virtual localization formula has become one of the main tools in enumerative geometry.

Another more recent perspective, and the perspective of this paper, is to consider the Bia{\l}ynicki-Birula stratification, the stratification of $X$ by the attracting manifolds $S^+_{\alpha}:= \{x \in X | \lim_{t \to 0} t \cdot x \in Z_\alpha\}$. Regarded as a structure on the stack $X/\bbC^\ast$, this is a very special case of a $\Theta$-stratification \cite[Def.~2.1.2]{Halpern-Leistner2014}, whose definition we recall in \Cref{S:main_result}. (It is also referred to as a KN stratification in the context of quotient stacks $X/G$ in \cite{Teleman2000TheRevisited}.) Given a $\Theta$-stratification of an algebraic stack $\mcX = \bigcup_\alpha \mcS_\alpha$, each stratum $\mcS_\alpha$ retracts onto a ``center" $\mcZ_\alpha$, which in the case of the Bia{\l}ynicki-Birula stratification is $Z_\alpha / \bbC^\ast$. This leads to a generalization of \ref{I:localization_homology} that is identical to \eqref{E:first_localization_formula}, except that $X$ is replaced with $\mcX$ and $Z_\alpha$ is replaced with $\mcZ_\alpha$. This formula was proved by Teleman-Woodward \cite[Eq.~1.12]{Teleman2009TheCurve} in the setting of a smooth quotient stack $X/G$, and it is closely related to the Jeffrey-Kirwan cohomological non-abelian localization formula \cite{Jeffrey1994LocalizationActions}.

The main contributions of this paper are:
\begin{itemize}
    \item We extend the $K$-theoretic non-abelian localization formula to any locally finitely presented algebraic derived stack with a $\Theta$-stratification in \eqref{E:perfect_formula}. This includes quasi-smooth stacks, which are the derived enhancements of stacks with a perfect obstruction theory. It also includes the more general class of derived stacks with perfect cotangent complexes, provided one can identify a ``fundamental cycle."

    \item We use this to prove a general virtual $K$-theoretic wall-crossing formula for quasi-smooth stacks that admit a proper good moduli space, \Cref{E:wall_crossing}, and we prove a new finiteness theorem for the cohomology of Atiyah-Bott complexes on the stack of one-dimensional sheaves on a smooth projective surface, \Cref{T:1D_sheaves}.

    \item We formulate and prove categorifications of the non-abelian localization theorem: We introduce a ``completion'' $\Perf(\mcX)_\beta^\wedge$ of the category of perfect complexes on $\mcX$ that satisfies an analogue of \ref{I:localization_cohomology}, \Cref{P:perf_completion}, and we introduce a category of ``highest weight cycles'' $\Cycles(\mcX)^{<\infty}$ that satisfies analogues of \ref{I:localization_modules} and \ref{I:localization_homology}, \Cref{T:non-abelian_localization}. Both are subcategories of the $\infty$-category $\QC(\mcX)$ of quasi-coherent complexes on $\mcX$. Our main theorem is a relative statement for a morphism $\rho : \mcX \to \mcB$, making it suitable for stacks, like the stack of Higgs bundles on a curve, that are defined over an affine base.
\end{itemize}

The non-abelian localization formula is more flexible than the original, because a stack can have many different $\Theta$-stratifications. For example, given an ample class $\ell \in \NS_{\bbC^\ast}(X)$, geometric invariant theory defines an equivariant stratification of $X$ consisting of $S_\alpha^+$ for every component $Z_\alpha$ such that $\ell|_{Z_\alpha}$ has negative weight, the repelling manifold $S^-_{\alpha}:= \{x \in X | \lim_{t \to \infty} t \cdot x \in Z_\alpha\}$ for any $\alpha$ where $\ell|_{Z_\alpha}$ has positive weight, and the semistable stratum $X^{\rm{ss}}$, which is the complement of all the others. The localization formula then says that for any $E \in \Perf(X/\bbC^\ast)$,
\begin{equation}\label{E:abelian_wall_crossing}
    \chi(X/\bbC^\ast, E) = \chi(X^{\ell-\rm{ss}}/\bbC^\ast, E) + \sum_\alpha \chi\left(Z_\alpha/\bbC^\ast, \frac{E|_{Z_\alpha}}{e(\bbN_{Z_\alpha/X})}\right),
\end{equation}
where $\chi(X/\bbC^\ast, E) = \sum_i (-1)^i \dim (H^i(X,E)^{\bbC^\ast})$. The left-hand-side is manifestly independent of $\ell$, so this gives a wall-crossing formula for how $\chi(X^{\ell-\rm{ss}}/\bbC^\ast, E)$ varies with $\ell$. (We discuss the general wall-crossing formula in \Cref{S:wall_crossing}.)

This wall-crossing example also highlights a subtlety in interpreting the non-abelian localization theorem: The $Z_\alpha$ terms in \eqref{E:abelian_wall_crossing} must depend on the stratification in order for the right-hand-side to be independent of $\ell$, so it is not sufficient to regard the terms as elements of localized $K$-theory. If $q \in K^0_{\bbC^\ast}(\rm{pt})$ corresponds to the character of weight $-1$, then the correct interpretation is to take $\chi\left(Z_\alpha, \frac{E|_{Z_\alpha}}{e(\bbN_{Z_\alpha/X})}\right) \in \bbZ[q^{\pm 1}]_{1-q^n}$, expand it as a Laurent series in $q$ if the weight of $\ell|_{Z_\alpha}$ is negative, and as a Laurent series in $q^{-1}$ otherwise, and then to keep the constant, i.e., $q^0$, term. In general, $e(\bbN_{\mcZ_\alpha/\mcX})^{-1}$ is interpreted as an explicit quasi-coherent sheaf that depends on the splitting $\bbN_{\mcZ_\alpha/\mcX} \cong \bbN_{\mcZ_\alpha/\mcX}^+ \oplus \bbN_{\mcZ_\alpha/\mcX}^-$ into positive and negative weight pieces -- see \cite[\S 10]{Woodward2010MomentTheory} for an expository account. Our approach provides a more conceptual explanation of $e(\bbN_{\mcZ_\alpha/\mcX})^{-1}$, discussed in \Cref{S:euler_class}.

A closely related issue arises when categorifying the non-abelian localization theorem. We do not know of an operation on categories that corresponds to inverting an element of $K_0$. Instead, both $\Perf(\mcX)^\wedge_\beta$ and $\Cycles(\mcX)^{<\infty}$ are closer in spirit to completion than to localization. For example, when $\mcX = \bbC^n / \bbC^\ast$ with the usual scaling action, regarded as a single $\Theta$-stratum, then $\Perf(\mcX)^\wedge_\beta = \Cycles(\mcX)^{<\infty}$ is equivalent to the derived category of complexes of graded $\bbC[x_1,\ldots,x_n]$-modules, where $x_i$ has weight $-1$, whose homology is finite dimensional in every weight and vanishes in sufficiently high weight. This category has $K_0 \cong \bbZ(\!(q)\!)$, where $q$ corresponds to the free module generated in weight $-1$.

\begin{rem}
    For smooth stacks, one can use the main structure theorems of \cite{Halpern-Leistner2015TheQuotient, Halpern-Leistner2021DerivedConjecture} to construct an isomorphism $K^\ast(\Perf(\mcX)) \cong \bigoplus_{\alpha} K^\ast(\Perf(\mcZ_\alpha))$ without localization. However, this isomorphism does not respect any multiplicative structure -- it is not a homomorphism of rings or of $K^0(\Perf(\mcX))$-modules, nor is it a homomorphism of $K^0_G(\rm{pt})$-modules when $\mcX = X/G$.
\end{rem}

\begin{rem}
    A more direct categorification of the Atiyah-Segal localization theorem (as opposed to non-abelian localization) for smooth quotient stacks $X/G$ was established in \cite{Chen2020EquivariantSpaces}. It shows that if $z \in G$ is semisimple and $X^z$ denotes the \emph{classical} fixed locus, then the canonical morphism of derived loop stacks $\mcL(X^z/G^z) \to \mcL(X/G)$ becomes an isomorphism of derived stacks after taking the formal completion of both stacks along either the fiber of the residual gerbe of $z \in G/_{\rm{adj}} G$ or the fiber of its image $[z] \in G/\!/G$ under the canonical morphisms $\mcL(X/G) \to G/_{\rm{adj}}G \to G/\!/G$. This implies equivalences of any category of sheaves one would like on these derived stacks.
\end{rem}

\subsection{Categorical analogue of K-homology}

Given a morphism of noetherian algebraic (derived) stacks $\rho : \mcX \to \mcB$, we will consider the following thick stable full subcategory of quasi-coherent complexes, which we call \emph{$K$-homology cycles},
\[
	\Cycles(\mcX/\mcB) := \left\{ K \in \QC(\mcX) \left|\; \forall E \in \Perf(\mcX),\; \rho_\ast(K\otimes E) \in \DCoh(\mcB) \right. \right \}.
\]

If $f : \mcX \to \mcY$ is a quasi-compact quasi-separated (qc.qs.) morphism over $\mcB$, universally of finite cohomological dimension, then the base change formula $f_\ast(E) \otimes F \cong f_\ast(E \otimes f^\ast(F))$ (see \cite[App.~A]{MR4560539}) implies that $f_\ast : \QC(\mcX) \to \QC(\mcY)$ maps $\Cycles(\mcX/\mcB) \to \Cycles(\mcY/\mcB)$. If $f$ is in addition cohomologically proper, meaning $H^0(f_\ast(E)) \in \DCoh(\mcY)$ for all $E \in \DCoh(\mcX)$, and of finite Tor amplitude, then $f_\ast : \QC(\mcX) \to \QC(\mcY)$ preserves perfect complexes. In this case the projection formula also implies that $f^\ast$ maps $\Cycles(\mcY/\mcB)$ to $\Cycles(\mcX/\mcB)$. It is natural to regard $\mcX \mapsto \Perf(\mcX)$ as a categorified cohomology theory of $\mcB$-stacks, in which case $\mcX \mapsto \Cycles(\mcX/\mcB)$ can be regarded as a categorified homology theory for $\mcB$-stacks.

\begin{ex}
	If $\mcB$ is perfect in the sense of \cite{Ben-Zvi2010}, such as a quotient stack in characteristic $0$ or a qc.qs. scheme, and $\rho$ is of finite presentation, separated, and representable by algebraic spaces, then $\Cycles(\mcX/\mcB)$ is the $\infty$-category of complexes of coherent sheaves on $\mcX$ whose support is proper over $\mcB$ \cite[Thm.~3.0.2]{Ben-Zvi2017}.
\end{ex}

$\Cycles(\mcX/\mcB)$ is a much richer object for stacks with positive dimensional stabilizers. An object $K \in \DCoh(\mcX)$ lies in $\Cycles(\mcX/\mcB)$ if its support is cohomologically proper and of finite cohomological dimension over $\mcB$. However, $\Cycles(\mcX/\mcB)$ can be larger than $\DCoh$:

\begin{ex}
	If $\mcX = BG$ for a linearly reductive $k$-group $G$ for a field $k$, then $\Cycles(BG)$ is the category of representations for which every isotypical summand has finite dimension.
\end{ex}

\begin{ex}
	If $\mcX = \bbA_k^1 / \bbG_m$ and $\mcB = \Spec(k)$ for a field $k$, then $\QC(\mcX)$ is equivalent, via the Rees construction, to the category of diagrams $\cdots \to F_{i+1} \to F_i \to \cdots$. $R\Gamma(\mcX,-)$ is the functor that assigns such a diagram to $F_0$, any perfect complex is a sum of shifts of line bundles, and tensoring by a line bundle corresponds to shifting the indexing of the diagram. Therefore, $\Cycles(\mcX/\mcB) \subset \QC(\mcX)$ is the full subcategory of diagrams such that $F_i \in \Perf(k), \forall i$. The objects $\mcO_{\bbA^1 \setminus 0}$ and the local cohomology complex $R\Gamma_0 \mcO_{\bbA^1}$ lie in $\Cycles(\mcX/\mcB)$ but not in $\DCoh(\mcX)$.
\end{ex}

\begin{rem}
    The assignment $K \mapsto \rho_\ast(K \otimes (-))$ gives a tautological functor from $\Cycles(\mcX/\mcB)$ to the $\infty$-category of exact functors of $\Perf(\mcB)^\otimes$-module categories
    \begin{equation} \label{E:dual_is_functor_category}
		\Cycles(\mcX/\mcB) \to \Fun^{\exact}_{\Perf(\mcB)^\otimes}(\Perf(\mcX),\DCoh(\mcB)).
	\end{equation}
    If $\mcB$ is regular and $\mcX$ is a perfect stack, then one can use the equivalence \[\QC(\mcX) \cong \Fun^L_{\QC(\mcB)^\otimes}(\QC(\mcX),\QC(\mcB))\] to show that \eqref{E:dual_is_functor_category} is an equivalence. In a sense taking the right-hand-side of \eqref{E:dual_is_functor_category} as the definition of $\Cycles(\mcX/\mcB)$ is a more natural choice -- for instance, covariant functoriality and the projection formula would be formal consequences of this definition -- but it will be convenient for us to work with the more concrete definition as a subcategory of $\QC(\mcX)$.
\end{rem}

\subsection{Statements of main results}

Our main goal will be to understand how the structure of a $\Theta$-stratification $\mcX = \bigcup_\alpha \mcS_\alpha$ relative to $\mcB$ is reflected in $\Cycles(\mcX/\mcB)$. We will fix our technical hypotheses, for reference throughout the paper:

\begin{hypothesis}\label{H:standard}
We consider a commutative diagram of algebraic derived $1$-stacks
\[\xymatrix@R=5pt{\mcX \ar[rr]^\rho \ar[dr]_{\pi_\mcX} & & \mcB \ar[dl]^{\pi_\mcB} \\ & \mcR & },\]
where all morphisms are quasi-separated and locally almost of finite presentation, $\mcR$ is noetherian, and $\mcB$ is quasi-compact. We further assume that the relative inertia morphism $\mcX \times_{\mcX \times_{\mcR} \mcX} \mcX \to \mcX$ is separated with affine fibers, and likewise for $\mcB$. This is automatic if, for instance, the relative diagonal $\mcX \to \mcX \times_\mcR \mcX$ or $\mcB \to \mcB \times_{\mcR} \mcB$ is affine.
\end{hypothesis}

\begin{hypothesis}\label{H:lfp}
	In addition to \Cref{H:standard}, $\pi_\mcX$ is locally of finite presentation, and any point of $\mcX$ has a quasi-compact open neighborhood $\mcU \subset \mcX$ such that $\rho|_{\mcU} : \mcU \to \mcB$ is universally of finite cohomological dimension \cite[Def.~A.1.4]{MR4560539}.
\end{hypothesis}

The finite presentation hypothesis is stronger in the derived setting than in the classical setting. A morphism of algebraic derived stacks is locally of finite presentation if and only if the underlying morphism of classical stacks is locally of finite presentation in the classical sense and $\bbL_{\mcX/\mcR} \in \Perf(\mcX)$ \cite[Thm.~7.4.3.18]{LurieHA}. On the other hand, if $f: \mcX \to \mcY$ is a morphism of locally noetherian classical stacks that is locally of finite presentation, then $f$ is locally \emph{almost} of finite presentation when regarded as a morphism of derived stacks. The finite cohomological dimension condition, which means that for some $d$, $\rho_\ast[d]$ is right $t$-exact after arbitrary base change, holds automatically if $\mcR$ is defined over $\bbQ$ \cite[Thm.~1.4.2]{Drinfeld2013OnStacks}. It implies that $\rho_\ast : \QC(\mcU) \to \QC(\mcB)$ commutes with filtered colimits \cite[Prop.A.1.5]{MR4560539}.

Our first main result assumes \Cref{H:lfp} and additionally that $\mcX$ is quasi-compact. Every object $E \in \QC(\mcZ_\alpha)$ decomposes functorially as a direct sum $\bigoplus_{w \in \bbZ} E^w$. The subcategory of highest weight cycles $\Cycles(\mcZ_\alpha/\mcB)^{<\infty} \subset \QC(\mcZ_\alpha)$ is the full subcategory of objects $E \cong \bigoplus_w E^w \in \Cycles(\mcZ_\alpha/\mcB)$ such that $E^w \cong 0$ for $w \gg 0$. In \Cref{D:highest_weight_cycles_multiple} we introduce a category of highest weight complexes $\QC(\mcX)^{<\infty}$ and highest weight cycles $\Cycles(\mcX/\mcB)^{<\infty} \subset \QC(\mcX)^{<\infty}$. An object $F\in \QC(\mcX)^{<\infty}$ lies in $\Cycles(\mcX/\mcB)^{<\infty}$ if its $!$-restriction to every stratum $\mcS_{\alpha}$ followed by $\ast$-restriction to $\mcZ_\alpha$ lies in $\Cycles(\mcZ_\alpha/\mcB)$.

We have canonical morphisms $\tot_\alpha : \mcZ_\alpha \to \mcX$ for each $\alpha$. We consider the pushforward functors $(\tot_\alpha)_\ast : \QC(\mcZ_\alpha) \to \QC(\mcX)$, and introduce in \Cref{D:sharp_pullback} a modified pullback functor $\tot^\sharp : \QC(\mcX) \to \QC(\mcZ_\alpha)$, which up to tensoring with a line bundle is the composition of $!$-pullback to $\mcS_\alpha$ and $\ast$-pullback to $\mcZ_\alpha$. If $P \in \Perf(\mcX)$, then $\tot^\sharp((-) \otimes P) \cong \tot^\sharp(-) \otimes \tot^\ast(P)$. Our main theorem categorifying and generalizing \ref{I:localization_modules} and \ref{I:localization_homology} is the following:
\begin{thm}[\Cref{T:non-abelian_localization}] \label{T:intro_main}
In the context of \Cref{H:lfp}, suppose $\mcX$ is quasi-compact and equipped with a $\Theta$-stratification relative to $\mcB$, $\mcX = \mcS_0 \cup \cdots \cup \mcS_N$. (See \Cref{D:relative_theta_stratification}.) Then $\tot_\ast$ and $\tot^\sharp$ preserve categories of highest weight cycles, and the resulting homomorphism on $K$-groups is an isomorphism
\[
\tot_\ast = \bigoplus_\alpha (\tot_\alpha)_\ast : \bigoplus_\alpha K_0(\Cycles(\mcZ_\alpha/\mcB)^{<\infty}) \xrightarrow{\cong} K_0(\Cycles(\mcX/\mcB)^{<\infty}).
\]
For any $[E] \in K_0(\Cycles(\mcX/\mcB)^{<\infty})$ we have $[E] = \sum_\alpha (\tot_\alpha)_\ast \left(\frac{[\tot^\sharp(E)]}{e(\bbN_{\mcZ_\alpha/\mcX})} \right)$.
\end{thm}

We say $\pi_\mcX$ is \emph{quasi-smooth} if the relative cotangent complex $\bbL_{\mcX/\mcR}$ is perfect of Tor-amplitude in $[-1,1]$. Because many stacks of interest, such as stacks of sheaves on smooth projective varieties of dimension $>2$, have cotangent complexes that are perfect of larger Tor amplitude, we were careful to pinpoint exactly where the quasi-smoothness is useful in the localization formula. As the theorem above shows, non-abelian localization holds as long as one can identify some class $\mathbf{O} \in \Cycles(\mcX/\mcB)^{<\infty}$ playing the role of the fundamental class, and one can compute $\tot^\sharp(\mathbf{O})$. If $\pi_\mcX$ is quasi-smooth then $\tot_\alpha^\sharp(\mcO_\mcX) \cong \mcO_{\mcZ_\alpha}$, and one can often take the structure sheaf $\mcO_\mcX$ as a fundamental class, giving the simpler expression in \Cref{T:quasi-smooth_localization}. In \Cref{S:wall_crossing}, we use this to give a universal wall-crossing formula for quasi-smooth stacks admitting a proper good moduli space. Although this is familiar to experts in wall-crossing for moduli of objects in abelian categories, we believe it is helpful to see a precise general formulation.

In \Cref{S:infinite_stratifications} we generalize the definition of $\Cycles(\mcX/\mcB)^{<\infty}$ to situations where $\mcX$ is not quasi-compact, in which case the stratification can be infinite. \Cref{P:nonabelian_localization_infinite_strata} is an analogue of \ref{I:localization_homology} in this setting. In \Cref{S:1D_sheaves} we use this to show that on the stack of pure sheaves of dimension $1$ on a smooth projective surface, certain ``admissible'' complexes have finite dimensional cohomology. This was known for closely related stacks, such as the stack of Higgs bundles on a smooth and proper curve \cite{halpernleistner2016equivariantverlindeformulamoduli}, and in those cases there are interesting Verlinde-type formulas for the virtual dimension of the space of sections. It would be very interesting if a similar formula holds for the stack of pure 1D-sheaves on a surface.

Finally, \Cref{T:main_relative} gives conditions on a morphism $f : \mcX \to \mcY$ between two stacks with $\Theta$-stratifications relative to $\mcB$ that imply that $f_\ast : \QC(\mcX) \to \QC(\mcY)$ preserves the categories $\Cycles(-/\mcB)^{<\infty}$. This is the non-abelian analogue of the relative form of the trace formula \cite[Thm.~3.5]{Thomason1992UneAlgebrique}.

\subsection{Notation}

By derived stacks, we mean functors from the $\infty$-category of simplicial commutative rings to $\infty$-groupoids that satisfy derived \'etale descent. We consider only $1$-stacks, meaning $\pi_i(\mcX(R)) \cong 0$ for $i>1$ for discrete simplicial commutative rings $R$, i.e., classical rings, in which case the restriction of $\mcX$ to the full subcategory of classical rings is a classical stack. Our stacks will mostly be algebraic, which under mild hypotheses is equivalent to saying they admit a derived cotangent complex and their restriction to classical rings is a classical algebraic stack.

For an algebraic derived stack $\mcX$, $\QC(\mcX)$ will denote the $\infty$-category obtained by Kan-extending the functor $R \mapsto R\Mod$ from the subcategory of affine derived stacks to all derived stacks. $\QC(\mcX)$ is a symmetric monoidal category, and $R\uHom_\mcX(E,F)$ denotes the internal Hom, i.e., it is characterized by the functorial isomorphism $\Map(G,R\uHom(E,F)) \cong \Map(G \otimes E,F)$.

\subsection{Author's note}

This paper is a strengthened version of two of my unpublished preprints, \cite{categorification_note} and \cite{halpern2015remarks}, and this represents the final finished form of those results. The quasi-smooth non-abelian localization theorem was already used in \cite{halpernleistner2016equivariantverlindeformulamoduli, halpernleistner2023structuremoduligaugedmaps}, so in addition to future applications, this paper puts those papers on firmer footing.

I would like to thank Harrison Chen, Martijn Kool, Davesh Maulik, Andrei Okounkov for inspiring discussions on these topics over the years. In addition, I would like to thank a group of talented young researchers for motivating me to complete this paper with their recent interest in non-abelian localization: Chenjing Bu, Ivan Karpov, Tasuki Kinjo, Henry Liu, Miguel Moreira, and Andr\'{e}s Ib\'{a}\~nez-N\'{u}\~nez. I would especially like to thank Miguel Moreira for catching an important error in the first version of this paper.

This work was supported by a Sloan Foundation Research Fellowship, the NSF FRG grant DMS-2052936, and the NSF CAREER grant DMS-1945478.

\section{Baric completion of a \texorpdfstring{$\Theta$}{Theta}-stratum}

\subsection{A brief review of \texorpdfstring{$\Theta$}{Theta}-stratifications}

Let $\Theta := \bbA^1 / \bbG_m$, where the coordinate function $t$ on $\bbA^1$ has weight $-1$ for the $\bbG_m$-action. Given a stack $\mcX$ as in \Cref{H:standard}, the stacks $\Filt(\mcX) := \Map(\Theta,\mcX)$ and $\Grad(\mcX):= \Map(B\bbG_m,\mcX)$ are also stacks satisfying \Cref{H:standard} by \cite[Prop.~1.1.2]{Halpern-Leistner2014}, where all of these mapping stacks denote the mapping stacks relative to $\mcR$. We have canonical morphisms
\begin{equation}\label{E:theta_stratum_arrows}
	\xymatrix{\Grad(\mcX) \ar@/^20pt/[rr]^-{\tot} \ar@/^10pt/[r]_-{\spfilt} & \Filt(\mcX) \ar@/^10pt/[l]^-{\gr} \ar[r]_-{\ev_1} & \mcX}.
\end{equation}
Here $\gr$ is induced by composition with the inclusion $\{0\}/\bbG_m \hookrightarrow \Theta$, $\spfilt$ is induced by the projection $\Theta \to B\bbG_m$, $\ev_1$ is induced by restriction to the point $1 \in \bbA^1$, and $\tot := \ev_1 \circ \spfilt$.

\begin{defn}[Relative $\Theta$-stratum]\label{D:relative_theta_stratum}
	A $\Theta$-stratum \cite[Def.~2.1]{Halpern-Leistner2014} relative to $\rho : \mcX \to \mcB$ is an open and closed substack $\mcS \subset \Filt_\mcR(\mcX)$ such that $\ev_1 : \mcS \to \mcX$ is a closed immersion and for any $f \in \mcS$, the composition of group homomorphisms $\bbG_m \to \Aut_{\mcX}(f(0)) \to \Aut_{\mcB}(\rho(f(0)))$ is trivial. By \cite[Lem.~1.3.8]{Halpern-Leistner2014}, $\mcS = \gr^{-1}(\mcZ)$ for the open and closed substack $\mcZ := \spfilt^{-1}(\mcS) \subset \Grad(\mcX)$, and we call $\mcZ$ the \emph{center} of $\mcS$. All of the morphisms in \eqref{E:theta_stratum_arrows} preserve $\mcZ$ and $\mcS$.
\end{defn}

We note that $\mcZ \subset \Grad(\mcX)$ is classified by a morphism $\mcZ \times B\bbG_m \to \mcX$, which lifts to both $\mcS$ and $\mcZ$ and also defines open and closed immersion $\mcZ \subset \Filt(\mcZ)$ and $\mcZ \subset \Filt(\mcS)$. Giving a morphism $\mcZ \times B\bbG_m \to \mcX$ is also equivalent to giving a morphism $\mcZ \to \mcX$ and a homomorphism of group schemes $(\bbG_m)_{\mcZ} \to I_{\mcX}|_{\mcZ}$ over $\mcZ$. Any quasi-coherent complex $E \in \QC(\mcZ)$ decomposes uniquely into a direct sum $\bigoplus_{w \in \bbZ} E^w$ such that $\bbG_m$ acts with weight $w$ on $E^w$.

The condition on $\Aut_{\mcB}(\rho(f(0)))$ ensures that under the morphism $\Filt(\mcX) \to \Filt(\mcB)$ induced by $\rho$, $\mcS$ lies over the image of the open and closed substack parameterizing trivial filtrations in $\mcB$, which can be canonically identified with $\mcB$ via the forgetful morphism $\ev_1$. It also guarantees that for any morphism $\mcB' \to \mcB$, $\mcS$ induces a $\Theta$-stratum in $\mcX' := \mcX \times_{\mcB} \mcB'$ in the sense of \cite[Def.~2.3.1]{Halpern-Leistner2014}, i.e., if $\mcS \times_{\mcB} \mcB'$ denotes the fiber product with respect to $\rho \circ \ev_1$, then $\mcS \times_{\mcB} \mcB' \to \Filt(\mcX')$ has an open and closed image that is a $\Theta$-stratum in $\mcX'$ relative to $\mcB'$ \cite[Cor.~1.3.16]{Halpern-Leistner2014}.

The interpretation as a mapping stack equips $\mcS$ canonically with the structure of a derived algebraic stack, which need not be classical even when $\mcX$ is. Many of the results we discuss make use of this derived structure, and particularly the following consequence for cotangent complexes, which is {\cite[Lem.~1.3.2 \& Lem.~1.5.5]{Halpern-Leistner2021DerivedConjecture}}:
\begin{lem}\label{L:stratum_cotangent}
	$\bbL_{\mcZ/\mcS} \cong \bigoplus_{w<0} \left(\tot^\ast(\bbL_{\mcX/\mcR})\right)^w [1]$ and $\spfilt^\ast(\bbL_{\mcS/\mcX}) \cong \bigoplus_{w>0} \left(\tot^\ast(\bbL_{\mcX/\mcR})\right)^w [1]$.
\end{lem}
Because we are led to derived algebraic geometry either way, we have allowed $\mcX$ to be a derived stack from the beginning.

\begin{defn}[Relative $\Theta$-stratification]\label{D:relative_theta_stratification}
A \emph{$\Theta$-stratification} of $\mcX$ is a totally preordered collection of disjoint open substacks $\bigsqcup_{i \in I} \mcS_
i \subset \Filt_\mcR(\mcX)$ such that for all $\alpha$:
\begin{enumerate}
    \item $\bigcup_{j>i} \ev_1(|\mcS_j|) \subset |\mcX|$ is closed. We call its open complement $\mcX_{\leq i} \subset \mcX$.
    \item $\mcS_i$ lies in the open substack $\Filt(\mcX_{\leq i}) \cong (\tot \circ \gr)^{-1} ( \mcX_{\leq i}) \subset \Filt(\mcX)$, and it is a $\Theta$-stratum in $\mcX_{\leq i}$ relative to $\mcB$. (See \Cref{D:relative_theta_stratum}.)
    \item $\ev_1 : \bigsqcup_{i \in I} \mcS_i \to \mcX$ is surjective (hence universally bijective).
\end{enumerate}
We will let $\mcZ_i := \spfilt^{-1}(\mcS_i) \subset \Grad(\mcX_{\leq i}) \subset \Grad(\mcX)$ denote the center of each stratum. Sometimes, we will abuse notation and regard $\mcS_i$ as a closed substack of $\mcX_{\leq i}$ via the morphism $\ev_1$.
\end{defn}

The condition (3) is a notational convenience, and differs slightly from the conventions of \cite[Def.~2.1.2]{Halpern-Leistner2014}, which only requires universal injectivity of $\ev_1$ and defines a semistable locus $\mcX^{\rm ss} := \mcX \setminus \bigcup_{j} \ev_1(\mcS_j) \subset \mcX$. In this paper, if $\mcX^{\rm{ss}} \neq \emptyset$, we adjoin a formal minimal element $0 \in I$ and let $\mcS_0 \subset \Filt(\mcX)$ denote the open substack parameterizing trivial filtrations in $\mcX^{\rm{ss}}$. Then $\ev_1 : \mcS_0 \cong \mcX^{\rm{ss}}$ is an isomorphism, i.e., we may regard the semistable locus as an additional (trivial) $\Theta$-stratum.

\subsection{The baric structure on a \texorpdfstring{$\Theta$}{Theta}-stratum}\label{S:baric_stratum}

Consider a stable $\infty$-category $\mcC$. A \emph{baric structure} on $\mcC$ consists of a collection of full stable cocomplete subcategories $\mcC^{\geq w}$ indexed by $w \in \bbZ$, such that $\mcC^{\geq w} \subseteq \mcC^{\geq w-1}$ for all $w$, and the inclusion $\mcC^{\geq w} \hookrightarrow \mcC$ admits a right adjoint $\radj{w} : \mcC \to \mcC^{\geq w}$ \cite{Achar2011BaricSheaves}. In this case we let $\mcC^{<w}$ be the right orthogonal complement of $\mcC^{\geq w}$, and
\[
	\ladj{w}(F) = \cofib (\radj{w}(F) \to F)
\]
is left adjoint to the inclusion $\mcC^{< w} \hookrightarrow \mcC$. Another way to say this is that one has a semiorthogonal decomposition $\mcC = \langle \mcC^{<w}, \mcC^{\geq w} \rangle$ for all $w \in \bbZ$. We will denote $\mcC^w := \mcC^{\geq w} \cap \mcC^{<w+1}$, and let $\beta^w := \radj{w} \circ \ladj{w+1} \cong \ladj{w+1} \circ \radj{w}$ be the projection functor.

If $\mcC$ is presentable, then we say that the baric structure is \emph{continuous} if $\radj{w}$ preserves filtered colimits for any $w$. This is equivalent to $\mcC^{<w} \subset \mcC$ being closed under filtered colimits in $\mcC$. To see this, observe that in the exact triangle $\colim_{\alpha} \radj{w}(E_\alpha) \to \colim_\alpha E_\alpha \to \colim_\alpha \ladj{w}(E_\alpha)$, the third term lies in $\mcC^{<w}$ if and only if this triangle is isomorphic to the exact triangle $\radj{w}(\colim_\alpha E_\alpha) \to \colim_\alpha E_\alpha \to \ladj{w}(\colim_\alpha E_\alpha)$. (The dual conditions, that $\mcC^{\geq w}$ is closed under small colimits and $\ladj{w} : \mcC \to \mcC^{<w}$ preserves small colimits, hold automatically for any baric decomposition of $\mcC$.) \emph{All of the baric structures on presentable stable $\infty$-categories we discuss in this paper will be continuous.}

It is sometimes convenient to consider baric structures indexed by $w \in \bbR$, in which case we also require that for any $w$, $\mcC^{<v} = \mcC^{<w}$ for $v$ in some interval $(w-\epsilon,w]$.

Pre-composition with the addition map $\Theta \times \Theta \to \Theta$ defines a morphism $a : \Theta \times \Filt(\mcX) \to \Filt(\mcX)$. This defines an action of the monoid $\Theta$ in the homotopy category of derived $\mcB$-stacks, i.e., $a$ satisfies the associativity and identity axioms up to homotopy. This action preserves any $\Theta$-stratum $\mcS \subset \Filt(\mcX)$. We call this a weak $\Theta$-action, and in \cite[Prop.~1.1.2]{Halpern-Leistner2021DerivedConjecture} we construct a (continuous) baric structure on $\QC(\mcS)$ for any stack $\mcS$ with a weak $\Theta$-action. The right adjoint to the inclusion $\QC(\mcS)^{\geq w} \hookrightarrow \QC(\mcS)$ is given by
\[
	\radj{w}(F):= \pi_\ast (\mcO_\Theta\langle w \rangle \otimes a^\ast(F)),
\]
where $\pi : \Theta \times \mcS \to \mcS$ is the projection, and $\mcO_\Theta\langle w \rangle$ is the line bundle corresponding to the free graded $\bbZ[t]$-module $\bbZ[t] \cdot t^w \subset \bbZ[t^{\pm 1}]$. $\radj{w}$ preserves $\APerf(\mcS)$ and $\Perf(\mcS)$, and therefore it restricts to a baric structure on these subcategories as well.

Now consider a $\rho : \mcX \to \mcB$ as in \Cref{H:standard}, and let $\mcS \subset \Filt_\mcR(\mcX)$ be a $\Theta$-stratum relative to $\mcB$. In \cite[Prop.~1.7.2]{Halpern-Leistner2021DerivedConjecture}, we construct a baric structure on $\QC(\mcX)$ characterized as the unique \emph{continuous} baric structure such that
\begin{itemize}
    \item $(\ev_1)_\ast : \QC(\mcS) \to \QC(\mcX)$ commutes with baric truncation, and
    \item the truncation functors $\radj{w}_\mcS$ and $\ladj{w}$ have locally bounded below homological $t$-amplitude in the sense that for any morphism $p : \Spec(A) \to \mcX$, there is a $d$ such that $p^\ast \circ \radj{w} (\QC(\mcX)_{\geq 0}) \subset \QC(\Spec(A))_{\geq d}$ and likewise with $\ladj{w}$ instead of $\radj{w}$.
\end{itemize}
Objects in $\QC(\mcX)^{\geq w}$ are set-theoretically supported on $\mcS$, hence we denote it $\QC_\mcS(\mcX)^{\geq w}$. The baric structure is $\mcB$-linear in the sense that if $\rho : \mcX \to \mcB$ is the structure morphism and $E \in \QC(\mcB)$, then $\rho^\ast(E) \otimes \QC(\mcX)^{<w} \subset \QC(\mcX)^{<w}$ and $\rho^\ast(E) \otimes \QC_{\mcS}(\mcX)^{\geq w} \subset \QC_{\mcS}(\mcX)^{\geq w}$ for any $w \in \bbZ$. It follows from \cite[Prop.~1.7.2(5)]{Halpern-Leistner2021DerivedConjecture} that the baric structure is local over $\mcB$ in the sense that if $\{\mcU_\alpha \to \mcB\}_{\alpha \in I}$ is a flat cover of $\mcB$, then $E \in \QC(\mcX)$ lies in $\QC_\mcS(\mcX)^{\geq w}$ or $\QC(\mcX)^{<w}$ if and only if its restriction to each $\mcU_\alpha \times_{\mcB} \mcX$ lies in the corresponding subcategory for the morphism $\mcU_\alpha \times_{\mcB} \mcX \to \mcU_\alpha$.

\begin{lem}\label{L:constant_weight_category}
	The pullback functor $\gr^\ast : \QC(\mcZ)^w \to \QC(\mcS)^w$ and the pushforward $(\ev_1)_\ast : \QC(\mcS)^w \to \QC(\mcX)^w$ are equivalences of $\infty$-categories.
\end{lem}
\begin{proof}
	This is proved in \cite[Thm.~2.2.3(3)]{Halpern-Leistner2021DerivedConjecture} for the categories $\APerf$, but the proof applies verbatim to $\QC$.
\end{proof}

\begin{lem}\label{L:first_local_cohomology_formula}
	$\forall F \in \QC(\mcX)$, the canonical maps $\radj{w}_\mcS(F) \to F$ determine an isomorphism
	\[
		\underset{w \to -\infty}{\colim} \radj{w}_\mcS(F) \xrightarrow{\cong} R\Gamma_\mcS(F).
	\]
\end{lem}
\begin{proof}
	The claim is unchanged if we replace $F$ with $R\Gamma_\mcS(F)$, so we may assume $F \in \QC_\mcS(\mcX)$. Then
	\[
		\cofib(\underset{w \to -\infty}{\colim} \radj{w}_\mcS(F) \to F) \in \bigcap_{w \in \bbZ} \QC_\mcS(\mcX)^{<w},
	\]
	so it suffices to show that the latter category is $0$. By \cite[Prop.~1.7.2(6)]{Halpern-Leistner2021DerivedConjecture}, an object $F \in \QC_\mcS(\mcX)$ lies in $\QC_\mcS(\mcX)^{<w}$ if and only if the same is true for all $H_n(F)$, so it suffices to show that $\QC(\mcX)_{<\infty} \cap \bigcap_{w \in \bbZ}(\QC_\mcS(\mcX)^{<w}) = 0$.

	To show that this category vanishes, it suffices to prove the claim of the lemma for all $F \in \QC_\mcS(\mcX)_{<\infty}$. In this case the filtration of $\radj{w}_\mcS(F)$ by $n^{th}$ derived infinitesimal neighborhoods in \cite[Eq.~11]{Halpern-Leistner2021DerivedConjecture} (see also \Cref{L:local_cohomology_formula} below) reduces the claim to objects of the form
	\[
		R\uHom_{\mcX}((\ev_1)_\ast(\Sym^n(\bbL_{\mcS/\mcX})),F) \cong (\ev_1)_\ast \left( R\uHom_{\mcS}(\Sym^n(\bbL_{\mcS/\mcX}), \ev_1^!(F)) \right),
	\]
	where $\bbL_{\mcS/\mcX}$ is the relative cotangent complex of the inclusion $\mcS \hookrightarrow \mcX$. Finally, we observe that the claim of the lemma holds for any object of the form $(\ev_1)_\ast(G)$ by \cite[Prop.~1.7.2(2), Prop.~1.1.2(2)]{Halpern-Leistner2021DerivedConjecture}.
\end{proof}

\begin{cor}\label{C:local_category}
    The pushforward functor $i_\ast : \QC(\mcX \setminus \mcS) \to \QC(\mcX)$ induces an equivalence $\QC(\mcX \setminus \mcS) \cong \bigcap_w \QC(\mcX)^{<w}$
\end{cor}
\begin{proof}
    $i_\ast$ is fully faithful, and for any $E \in \QC_{\mcS}(\mcX)^{\geq w}$, $\Hom(E,i_\ast(F)) \cong \Hom(i^\ast(E),F) \cong 0$ because $i^\ast(E) \cong 0$. Therefore $i_\ast(F) \in \QC(\mcX)^{<w}$ for any $w$. On the other hand, for any $F \in \bigcap_w \QC(\mcX)^{<w}$, \Cref{L:first_local_cohomology_formula} implies that $R\Gamma_\mcS(F) \cong \colim_w \radj{w}_\mcS(F) \cong 0$, so $F \in i_\ast(\QC(\mcX \setminus \mcS))$.
\end{proof}

We introduce the following convenient terminology:
\begin{defn}[Highest weight complexes]
For $F \in \QC(\mcX)$, we let
\begin{equation}\label{E:highest_weight}
\hwt(F) := \sup \left\{w \mid \exists d \in \bbZ \text{ s.t. } \radj{w}(\ev_1^!(\tau_{\leq d}(R\Gamma_{\mcS}(F)))) \neq 0 \right\} \in \bbZ \cup {\pm \infty},
\end{equation}
and for $F \in \QC(\mcS)$ or $\QC(\mcZ)$ we let $\hwt(F) := \sup \{w \mid \radj{w}(F) \neq 0\}$. We define $\QC(\mcX)^{<\infty} := \{F \in \QC(\mcX) \mid \hwt(F)<\infty\}$ and likewise for $\QC(\mcZ)^{<\infty}$ and $\QC(\mcS)^{<\infty}$.
\end{defn}

Recall that an algebraic derived stack $\mcX$ is eventually coconnective if $H_i(\mcO_\mcX) \cong 0$ for all $i \gg 0$.
\begin{lem}\label{L:highest_weight_condition}
For $F \in \QC(\mcX)$ and $v \in \bbZ \cup \{\infty\}$, $F \in \QC(\mcX)^{<v}$ if and only if $\hwt(F) < v$, and the analogous claim holds for $\mcS$ and $\mcZ$. If $\mcX$ is eventually coconnective or $F \in \QC(\mcX)_{<\infty}$, then $\hwt(F) = \hwt(\ev_1^!(F))$.
\end{lem}
\begin{proof}
    This is explained in \cite[Rem.~1.7.3 \& \S 1.7.1]{Halpern-Leistner2021DerivedConjecture}. The reason for the homological truncation $\tau_{\leq d}$ in \eqref{E:highest_weight} is that $(\ev_1)_\ast(\QC(\mcS)^{\geq w})$ does not generate $\QC_\mcS(\mcX)^{\geq w}$ under extensions, shifts, and filtered colimits alone. When $\mcX$ is not eventually coconnective, it is also necessary to take limits of towers $\to \cdots \to F_1 \to F_0$ such that $\tau_{\leq k}(F_i)$ is eventually constant in $i$ for any $k$.
\end{proof}

In \cite[Appendix A]{Halpern-Leistner2021DerivedConjecture}, for any closed substack $i : \mcS \hookrightarrow \mcX$, we construct an ascending chain of closed substacks
\[
\mcS \hookrightarrow \mcS^{(1)} \hookrightarrow \mcS^{(2)} \hookrightarrow \cdots \hookrightarrow \mcX
\]
such that each $\mcS^{(n)} \hookrightarrow \mcS^{(n+1)}$ is surjective, each $\mcO_{\mcS^{(n)}} \in \APerf(\mcX)$, and there is a canonical exact triangle for all $n>0$,
\begin{equation}\label{E:infinitesimal_neighborhood}
	i_\ast(\Sym^n_\mcS(\bbL_{\mcS/\mcX}[-1])) \to \mcO_{\mcS^{(n)}} \to \mcO_{\mcS^{(n-1)}} \to.
\end{equation}
$\mcS^{(n)}$ is the analog in derived algebraic geometry of the $n^{th}$ infinitesimal neighborhood of $\mcS$ in $\mcX$. For $E \in \QC(\mcX)_{<\infty}$, \cite[Thm.~A.0.1]{Halpern-Leistner2021DerivedConjecture} implies the canonical map in $\QC(\mcX)$ is an isomorphism
\begin{equation}\label{E:local_cohomology_formula}
	\colim_{n \to \infty} R\uHom_\mcX(\mcO_{\mcS^{(n)}},E) \xrightarrow{\cong} R\Gamma_\mcS(E).
\end{equation}
This formula does not hold for unbounded complexes in general, but we have the following:

\begin{lem}\label{L:local_cohomology_formula}
	For any $E \in \QC(\mcX)^{<\infty}$, \eqref{E:local_cohomology_formula} is an isomorphism.
\end{lem}
\begin{proof}
	We first claim that $R\uHom_{\mcX}(\mcO_{\mcS^{(n)}},-)$ preserves $\QC(\mcX)^{<w}$ for any $n\geq 0$ and $w \in \bbZ$. By induction on \eqref{E:infinitesimal_neighborhood}, it suffices to show this for the functor \[R\uHom_{\mcX}((\ev_1)_\ast \Sym^n_\mcS(\bbL_{\mcS/\mcX}[-1]),-) \cong (\ev_1)_\ast R\uHom_{\mcS}(\Sym^n_\mcS(\bbL_{\mcS/\mcX}[-1]), \ev_1^!(-)).\]
	Because $(\ev_1)_\ast$ commutes with baric truncation, it suffices to show that $\ev_1^!$ maps $\QC(\mcX)^{<w}$ to $\QC(\mcS)^{<w}$, and that $R\uHom_{\mcS}(\Sym^n_\mcS(\bbL_{\mcS/\mcX}[-1]),-)$ maps $\QC(\mcS)^{<w}$ to $\QC(\mcS)^{<w-n}$. The first claim follows from the fact that the left adjoint $(\ev_1)_\ast$ of $\ev_1^!$ preserves $\QC(-)^{\geq w}$, and the second claim follows from the fact that $\Sym^n_\mcS(\bbL_{\mcS/\mcX}[-1]) \in \QC(\mcS)^{\geq n}$, by \Cref{L:stratum_cotangent}, and so $\Sym^n_\mcS(\bbL_{\mcS/\mcX}[-1]) \otimes (-)$ maps $\QC(\mcS)^{\geq w}$ to $\QC(\mcS)^{\geq w+n}$ by \cite[Prop.~1.1.2(6)]{Halpern-Leistner2021DerivedConjecture}.

	To prove that \eqref{E:local_cohomology_formula} is an isomorphism, it suffices by \Cref{L:first_local_cohomology_formula} to show this map is an isomorphism after applying $\radj{w}_\mcS$ for any $w \in \bbZ$, so it suffices to show that for any $w \in \bbZ$,
	\begin{equation}\label{E:intermediate_local_cohomology_formula}
		\colim_{n \to \infty} \radj{w}_\mcS(R\uHom_\mcX(\mcO_{\mcS^{(n)}},E)) \xrightarrow{\cong} \radj{w}_\mcS(E).
	\end{equation}
	The fact that $R\uHom(\mcO_{\mcS^{(n)}},-)$ preserves $\QC(\mcX)^{<w}$ implies that replacing $E$ with $\radj{w}_\mcS(E)$ has no effect on the left-hand-side. $\radj{w}_\mcS(E)$ is set-theoretically supported on $\mcS$, so it suffices to show that \eqref{E:intermediate_local_cohomology_formula} is an isomorphism for any $E \in \QC_\mcS(\mcX)^{<v}:= \QC_\mcS(\mcX) \cap \QC(\mcX)^{<v}$ for some $v\in \bbZ$, so we will assume this about $E$ for the remainder of the proof.

    Now, the fact that $R\uHom_{\mcS}(\Sym^n_\mcS(\bbL_{\mcS/\mcX}[-1]),-)$ maps $\QC(\mcS)^{<v}$ to $\QC(\mcS)^{<v-n}$ implies that the colimit on the left hand side of \eqref{E:intermediate_local_cohomology_formula} is constant for $n\geq v-w$. In particular, it suffices to show that for $E \in \QC_\mcS(\mcX)^{<v}$, the homomorphism $\radj{w}_\mcS R\uHom(\mcO_{\mcS^{(v-w-1)}}, E) \to \radj{w}_\mcS(E)$ is an isomorphism. If we define $I^{(n)} = \fib(\mcO_\mcX \to \mcO_{\mcS^{(n)}})$, then this is equivalent to
	\[
		0 \cong \radj{w}_{\mcS}(R\uHom(I^{(v-w-1)}, E)), \text{ i.e., } R\uHom(I^{(v-w-1)}, E) \in \QC(\mcX)^{<w}.
	\]
	From \eqref{E:local_cohomology_formula} we know that this holds for $E \in \QC(\mcX)_{<\infty}$, and because $\QC_\mcS(\mcX)^{<v}$ is closed under truncation \cite[Prop.~1.7.2(6)]{Halpern-Leistner2021DerivedConjecture}, one has $R\uHom(I^{(v-w-1)}, \tau_{\leq d}(E)) \in \QC(\mcX)^{<w}$ for all $d$. Because $R\uHom(I^{(v-w-1)}, -)$ commutes with limits and $\QC(\mcX)^{<w}$ is closed under limits, we see that the condition holds for $E = \lim_{d \to \infty} \tau_{\leq d}(E)$. 
\end{proof}

\begin{lem}
For all $v \in \bbZ$, $\ev_1^! : \QC(\mcX)^{<v} \to \QC(\mcS)^{<v}$ commutes with filtered colimits, where $\ev_1^! := R\uHom(\mcO_\mcS,-)$ is right adjoint to the pushforward functor $(\ev_1)_\ast : \QC(\mcS)^{<v} \to \QC(\mcX)^{<v}$.
\end{lem}
\begin{proof}
The right adjoint of $\ev_1^\ast : \QC(\mcS)^{\geq w} \to \QC(\mcX)^{\geq w}$ is $\radj{w} \circ \ev_1^!$. Because $(\ev_1)_\ast$ commutes with the inclusions $\QC(-)^{\geq w} \subseteq \QC(-)$, passing to right adjoints shows that
\[\radj{w} \circ \ev_1^! \cong \radj{w} \circ \ev_1^! \circ \radj{w}.\]
Now given a filtered diagram $\{F_\alpha \in \QC(\mcX)^{<v}\}_{\alpha \in I}$, we wish to show that the canonical homomorphism is an isomorphism
\[
\colim_{\alpha} \ev_1^!(F_\alpha) \to \ev_1^!(\colim_\alpha F_\alpha).
\]
By \cite[Prop.~1.1.2(2)]{Halpern-Leistner2021DerivedConjecture}, it suffices to show this is an isomorphism after applying $\radj{w}$ for any $w \in \bbZ$. It therefore suffices to show that
\[
\colim_\alpha \radj{w} \ev_1^!(\radj{w}(F_\alpha)) \to \radj{w} \ev_1^!(\colim_\alpha \radj{w}(F_\alpha)).
\]
Using the fact that $\radj{w}(F_\alpha)$ has a functorial finite filtration by $\radj{w'}(F_\alpha)$ for $w \leq w'<v$, it suffices to prove the claim for the associated graded pieces. We may therefore assume $v=w+1$, and $F_\alpha \in \QC(\mcX)^w$ for all $\alpha$. In this case, the functor $\radj{w} \ev_1^! : \QC(\mcX)^w \to \QC(\mcS)^w$ is the inverse of the isomorphism $(\ev_1)_\ast : \QC(\mcS)^w \to \QC(\mcX)^w$ in \Cref{L:constant_weight_category}, so it commutes with filtered colimits.
\end{proof}

\subsection{A canonical filtration for derived self-intersections}

We will use the construction of derived deformation to the normal cone $\mcD_{\mcX/\mcY}$ for a morphism $f : \mcX \to \mcY$ that is locally almost of finite presentation. This is developed in the soon-to-be-released \cite{normal_cone}, extending previous work for non-derived stacks \cite{MR2877433,MR4748168}. The deformation to the normal cone is defined as $\mcD_{\mcX/\mcY}:= \Map_{\mcY \times \Theta}(\mcY \times B\bbG_m, \mcX \times \Theta)$, and by \cite[Thm.~C]{normal_cone} it is an algebraic derived $2$-stack, meaning the $\infty$-groupoid of $T$-points for any classical, i.e., discrete, affine derived scheme $T$ is $2$-connected. It is locally almost of finite presentation over $\mcY \times \Theta$. It is a $1$-stack, i.e., a stack in the sense of classical algebraic geometry, if $f : \mcX \to \mcY$ is representable by Deligne-Mumford stacks. The algebraicity depends crucially on the fact that the closed immersion $B\bbG_m \to \Theta$ has finite Tor-amplitude. If one wishes to use earlier algebraicity results (e.g., \cite{MR4560539}), one could also prove the results of this section using the derived deformation to the tangent cone $\Map_{\mcY \times \Theta}(\mcY \times \mcW, \mcX \times \Theta)$, where $\mcW := \Spec(\bbZ[x,y]/(xy)) / \bbG_m \to \Theta$ is a finite flat morphism.

Concretely, a $T$-point of $\mcD_{\mcX/\mcY}$ consists of a morphism $T \to \mcY$, an invertible sheaf with a section $s : \mcO_T \to L$, and a morphism $V(s) \to \mcX$ over $\mcY$, where $V(s)$ is the derived vanishing locus of $s$. There is a canonical diagram
\[
	\xymatrix{\mcX \times \Theta \ar[r]^{\tilde{f}} & \mcD_{\mcX/\mcY} \ar[r] & \mcY \times \Theta}.
\]
The fiber over $1 \in \Theta$ is the diagram $\mcX \to \mcY = \mcY$, where the first morphism is $f$. The fiber of $\mcD_{\mcX/\mcY}$ over $\{0\} / \bbG_m \hookrightarrow \Theta$ is the linear stack $\bbV_{\mcX}(\bbL_{\mcX/\mcY}[-1]) / \bbG_m$, whose $T$ points consist of a $T$-point $\xi : T \to \mcX$, an invertible sheaf $L$ on $T$, and a point of $\Map_T(\xi^\ast(\bbL_{\mcX/\mcY}[-1]), L^{-1})$. The fiber of the diagram above over $\{0\}/\bbG_m \hookrightarrow \Theta$ is $\mcX \to \bbV_\mcX(\bbL_{\mcX/\mcY}[-1]) \to \mcY$, where the first morphism is the inclusion of the $0$-section and the second morphism is the projection to $\mcX$ composed with $f$. 

\begin{prop}\label{P:push_pull_filtration_closed_immersion}
	Let $i : \mcS \to \mcX$ be a closed immersion of locally noetherian algebraic derived stacks. Then for $E \in \QC(\mcS)_{<\infty}$, there is a functorial convergent filtration $i^!(i_\ast(E)) \cong \colim (E_0 \to E_{-1} \to \cdots)$ with $\gr_n(E_\bullet) := \cofib(E_{n+1} \to E_n) \cong R\Hom_{\mcS}^\otimes(\Sym^{-n}(\bbL_{\mcS/\mcX}), E)$ for all $n \in \bbZ$. Furthermore, if $i = \ev_1 : \mcS \to \mcX$ is the inclusion of a $\Theta$-stratum and $E \in \QC(\mcS)^{<w}_{<\infty}$, then $\radj{v}(E_n) \to \radj{v}(E)$ is an isomorphism for all $n \leq v-w+1$, so the filtration is finite. If in addition $\bbL_{\mcS/\mcX}$ is perfect, then this finite filtration exists for any $E \in \QC(\mcS)^{<w}$, without homological bounds.
\end{prop}

\begin{proof}

We first observe that $\tilde{i} : \mcS \times \Theta \to \mcD_{\mcS/\mcX}$ is a closed immersion. Indeed, the formation of $\mcD_{\mcS/\mcX}$ is smooth local over $\mcX$, so it suffices to show this when $\mcX = \Spec(A)$ is affine. The underlying classical stack of $\mcD_{\mcS/\mcX}$ is $\Spec(R(I))/\bbG_m$, where $I \subseteq A$ is the ideal of definition of $\mcS \hookrightarrow \Spec(A)$ and $R(I):= \bigoplus_{n \in \bbZ} I^n \cdot t^{-n} \subset A[t^{\pm 1}]$, and $\tilde{i}$ is the closed immersion defined by the ideal in $R(I)$ generated by $\sum_{n > 0} I^n t^{-n}$.

Because $\tilde{i}$ is a closed immersion, we can consider the functor $\tilde{i}^! (\tilde{i}_\ast((-) \boxtimes \mcO_\Theta)) : \QC(\mcS) \to \QC(\mcS \times \Theta)$. The formation of $R\uHom(F,E)$ commutes with base change along morphisms of finite Tor-amplitude if $F$ is almost perfect and $E \in \QC_{<\infty}$ \cite[Lem.~2.4.3]{Halpern-Leistner2021DerivedConjecture}, so the fiber of $\tilde{i}^! (\tilde{i}_\ast(E \boxtimes \mcO_\Theta))$ over $1 \in \Theta$ is isomorphic to $i^!(i_\ast(E))$, and the fiber over $\{0\}/\bbG_m$ is isomorphic to $\bigoplus_{n} R\uHom(\Sym_n(\bbL_{\mcS/\mcX}),E)$, where the $n^{th}$ term has weight $n$. Under the Rees correspondence for $\QC(\mcS \times \Theta)$, this corresponds to a diagram $\cdots \to E_1 \to E_0 \to E_{-1} \to \cdots$ with colimit $i^! (i_\ast(E))$ and associated graded pieces as in the statement of the theorem.


What remains is to show that $E_i \cong 0$ for $i > 0$, i.e., that $\tilde{i}^! (\tilde{i}_\ast(E \boxtimes \mcO_\Theta)) \in \QC(\mcS \times \Theta)^{<1}$. This is equivalent, by adjunction, to showing that $\Hom(\tilde{i}^\ast(\tilde{i}_\ast(G)), E \boxtimes \Theta) \cong 0$ for all $G \in \QC(\mcS \times \Theta)^{\geq 1}$. Because $E \boxtimes \mcO_\Theta \in \QC(\mcS \times \Theta)^{<1}$, it suffices to show that $\tilde{i}^\ast(\tilde{i}_\ast(-))$ maps $\QC(\mcS \times \Theta)^{\geq 1}$ to itself. The claim is again smooth local, so we may assume $\mcX$ is affine and therefore that we can write any $G \in \QC(\mcS \times \Theta)^{\geq 1}$ as a filtered colimit $G \cong \colim_\alpha G_\alpha$ with $G_\alpha \in \Perf(\mcS \times \Theta)$. Then $G \cong \radj{1}(G) \cong \colim_\alpha \radj{1}(G_\alpha)$. Because $\radj{1}(G_\alpha) \in \APerf(\mcS \times \Theta)^{\geq 1}$, it suffices to show that $\tilde{i}^\ast \circ \tilde{i}_\ast$ preserves $\APerf(\mcS)^{\geq 1}$. By \cite[Lem.~1.5.4]{Halpern-Leistner2021DerivedConjecture}, any $G \in \APerf(\mcS \times \Theta)$ lies in $\APerf(\mcS \times \Theta)^{\geq 1}$ if and only if its restriction lies in $\APerf(\mcS \times B\bbG_m)^{\geq 1}$. The result now follows from the observation that the fiber of $\tilde{i}^\ast \circ \tilde{i}_\ast$ over $\{0\}/\bbG_m \hookrightarrow \Theta$ is $\Sym(\bbL_{\mcS/\mcX}) \otimes (-)$, which preserves $\APerf(\mcS \times B\bbG_m)^{\geq 1}$ because $\bbL_{\mcS/\mcX}$ has weight $1$.

Finally, suppose $i = \ev_1$ for a $\Theta$-stratum, and $E \in \QC(\mcX)^{<w}_{<\infty}$. By \Cref{L:stratum_cotangent} $\Sym^{-n}(\bbL_{\mcS/\mcX}) \in \APerf(\mcS)^{\geq -n}$, and therefore $R\uHom_\mcS(\Sym^{-n}(\bbL_{\mcS/\mcX}), E) \in \QC(\mcS)^{<w+n}$. We have
\[
\radj{v}(\gr_n(E_\bullet)) \cong \radj{v}(R\uHom_\mcS(\Sym^{-n}(\bbL_{\mcS/\mcX}),E)) \cong 0
\]
if $w+n\leq v$, so $\radj{v}(E_{n+1}) \to \radj{v}(E)$ is an isomorphism for $n$ in that range. If $\bbL_{\mcS/\mcX}$ is perfect, then each of the functors $R\uHom_\mcS(\Sym^n(\bbL_{\mcS/\mcX}), E) \cong \Sym^n(\bbL_{\mcS/\mcX}^\vee) \otimes E$ are right $t$-exact up to a shift, and therefore so are the functors $\Phi_n$ that take $E \in \QC(\mcS)^{<w}$ to the term $E_n$ in the filtration of $i^!(i_\ast(E))$. An $E \in \QC(\mcS)$ lies in $\QC(\mcS)^{<w}$ if and only if $H^i(E) \in \QC(\mcS)^{<w}$, so the functor $\Phi_n : \QC(\mcS)^{<w}_{<\infty} \to \QC(\mcS)^{<w}$ extends uniquely to a functor $\QC(\mcS)^{<w} \to \QC(\mcS)^{<w}$ that is right $t$-exact up to a shift, given by $E \mapsto \lim_{d \to \infty}(\tau_{\leq d}(E))$.

\end{proof}

\begin{thm}\label{T:push_pull_filtration_stratum}
	Suppose that $\spfilt : \mcZ \to \mcS$ is the center of a $\Theta$-stratum. Then $\forall F \in \QC(\mcZ)$ there is a functorial diagram whose formation commutes with filtered colimits
	\begin{equation}\label{E:push_pull_filtration}
		\cdots \to F_2 \to F_1 \to F_0 \cong \spfilt^\ast(\spfilt_\ast(F)),
	\end{equation}
	such that $\forall n$, $\gr_n(F^\bullet):= \cofib(F_{n+1} \to F_n) \cong \Sym^n(\bbL_{\mcZ/\mcS}) \otimes F$. Furthermore, if $F \in \QC(\mcZ)^{<\infty}$, then for any $v \in \bbZ$, $\radj{v}(F_i) \cong 0$ for all $i \gg 0$. In particular, \eqref{E:push_pull_filtration} gives a finite filtration of $\radj{v}(\spfilt^\ast(\spfilt_\ast(F)))$ whose associated graded pieces are $\radj{v}(\Sym^n(\bbL_{\mcZ/\mcS}) \otimes F)$.
\end{thm}

\begin{proof}
	We consider the $\tilde{F}:=\tilde{\spfilt}_\ast(\tilde{\spfilt}^\ast(F \boxtimes \mcO_\Theta)) \in \QC(\mcZ \times \Theta)$. Under the Rees equivalence, $\tilde{F}$ corresponds to a diagram $\cdots \to F_1 \to F_0 \to F_{-1} \to \cdots$ in $\QC(\mcZ)$. The restriction of $\tilde{F}$ to the fiber over $1 \in \Theta$, which under the Rees equivalence corresponds to $\colim F_i$, is isomorphic to $\spfilt^\ast(\spfilt_\ast(F))$ by the derived base change formula. Likewise, by the derived base change and projection formula,
	\[\tilde{F}|_{\{0\}/\bbG_m} \cong \Sym(\bbL_{\mcZ/\mcS}) \otimes F \in \QC(\mcZ \times B\bbG_m),\]
	where $F$ is concentrated in weight $0$ for the grading coming from the $B\bbG_m$ factor, $\bbL_{\mcZ/\mcS}$ is concentrated in weight $1$, and $\Sym(\bbL_{\mcZ/\mcS})$ is the pushforward of the structure sheaf along $\bbV_\mcZ(\bbL_{\mcZ/\mcS}) \to \mcZ$. This implies that $F_n \to F_{n-1}$ is an isomorphism for all $n \geq 0$, hence $F_0 \cong \colim F_n \cong \spfilt^\ast(\spfilt_\ast(F))$.
	
	What remains is to prove the vanishing claim that $\radj{w}(\tilde{F}_n) \cong 0$ for $n \gg 0$. We first show that $\mcD_{\mcZ/\mcS}$ is also a $\Theta$-stratum, and that $\mcZ \times \Theta \to \mcD_{\mcZ/\mcS}$ is $\Theta$-equivariant for a certain $\Theta$-action on the source. To simplify notation, we let $\mcD := \mcD_{\mcZ / \mcS}$.
	
	The fact that $\mcZ$ is the center of a $\Theta$-stratum means it is equipped with an open and closed immersion $\mcZ \to \Grad(\mcZ)$, corresponding to a cocharacter $\gamma : (\bbG_m)_{\mcZ} \to I_{\mcZ}$. Consider the morphism
    \[
    \alpha : \mcZ \times B\bbG_m \to \Grad(\mcD)
    \]
    that classifies the closed immersion $\mcZ \times B\bbG_m \hookrightarrow \mcD$ over $\{0\}/\bbG_m \hookrightarrow \Theta$ along with the homomorphism of group schemes $\tilde{\gamma} = (\gamma(t),t) : (\bbG_m)_{\mcZ \times B\bbG_m} \to I_{\mcZ \times B\bbG_m} \cong I_{\mcZ} \times \bbG_m$ that is $\gamma$ on the left factor and the identity on the right factor $\bbG_m$. $\alpha$ is a closed immersion because it factors as an open and closed immersion $\mcZ \times B\bbG_m \to \Grad(\mcZ \times B\bbG_m)$ corresponding to the cocharacter $\tilde{\gamma}$, followed by the closed immersion $\Grad(\mcZ \times B\bbG_m) \to \Grad(\mcD)$ \cite[Cor.~1.1.8]{Halpern-Leistner2014}.

	We now consider the pullback of the cotangent complex $\bbL_{\mcD/ (\mcS \times \Theta)}$ under the morphism $\tot \circ \alpha : \mcZ \times B\bbG_m \to \mcD$. From the description of $\mcD$ as a Weil restriction, this morphism corresponds to a morphism $\xi : V(s) \to \mcZ$ over $\mcS$, where $i : V(s) \to \mcZ \times B\bbG_m$ is the \emph{derived} vanishing locus of the zero section $s = 0 : \mcO_\mcZ \to \mcO_{\mcZ} \langle -1 \rangle$ of the twist of $\mcO_\mcZ$ by a character of $\bbG_m$ of weight $1$ for the $B\bbG_m$ factor. Then from the general description of the cotangent complex of a Weil restriction \cite[Prop.~5.1.10]{MR4560539}, we have
    \[
    (\tot \circ \alpha)^\ast(\bbL_{\mcD/\mcS\times \Theta}) \cong i_+(\xi^\ast(\bbL_{\mcZ/\mcS})),
    \]
    where $i_+ \cong i_\ast(-)\otimes \mcO_\mcZ \langle -1 \rangle [-1]$ is the left adjoint of the pullback $i^\ast$. (The formation of $i_+$ commutes with base change, so this formula for $i_+$ is obtained by base change from an explicit calculation for $i : \{0\}/\bbG_m \hookrightarrow \Theta$.) Moreover, $\xi$ factors as the composition $V(s) \xrightarrow{i} \mcZ \times B\bbG_m \to \mcZ$, where the second map is projection onto the left factor. The projection formula then gives
    \[
     (\tot \circ \alpha)^\ast(\bbL_{\mcD/\mcS\times \Theta}) \cong \mcO_\mcZ\langle -1 \rangle [-1] \otimes \bbL_{\mcZ/\mcS}.
    \]
    As a complex on $\mcZ \times B\bbG_m$, this has strictly negative weights with respect to $\gamma$ (coming from $\bbL_{\mcZ/\mcS}$, by \Cref{L:stratum_cotangent}) and weight $1$ with respect to the right factor of $B\bbG_m$. The map $\Grad(\mcZ) \to \Grad(\mcZ)$ raising a cocharacter to a positive power is an open and closed immersion \cite[Cor.~1.3.11]{Halpern-Leistner2014}, so after replacing $\gamma$ with $\gamma^2$, we may assume that $\bbL_{\mcZ/\mcS}$ has highest $\gamma$-weight $<-1$.
	
	The calculations above show that, $\tot^\ast(\bbL_{\mcD/\mcS\times \Theta})$ has strictly negative $\tilde{\gamma}$-weight. This implies:
	\begin{itemize}
	\item $\alpha$ is \'etale, hence an open and closed immersion: After pulling back under $\tot \circ \alpha$ and taking the $\tilde{\gamma}$-weight $0$ summand, the third term in the exact triangle $\bbL_{\mcS \times \Theta}|_{\mcD} \to \bbL_{\mcD} \to \bbL_{\mcD/\mcS \times \Theta} \to $ vanishes. This identifies $\alpha^\ast(\bbL_{\Grad(\mcD)})$, which is the $\tilde{\gamma}$-weight-$0$ summand of $(\tot \circ \alpha)^\ast(\bbL_{\mcD})$ by \Cref{L:stratum_cotangent}, with $\bbL_{\mcZ \times B\bbG_m} \cong \bbL_{\mcZ} \oplus \mcO_\mcZ[-1]$, the $\tilde{\gamma}$-weight-$0$ summand of $(\tot \circ \alpha)^\ast(\bbL_{\mcS \times \Theta}|_{\mcD}) \cong \spfilt^\ast(\bbL_{\mcS}) \oplus (\mcO_{\mcZ}\langle 1 \rangle \oplus \mcO_{\mcZ}[-1])$.

	\item If we define $\tilde{\mcS} \subset \Filt(\mcD)$ to be the open and closed preimage of $\mcZ \times B\bbG_m \subset \Grad(\mcD)$ under $\gr : \Filt(\mcD) \to \Grad(\mcD)$, then $\ev_1 : \tilde{\mcS} \to \mcD$ is \'etale: In general $\spfilt^\ast(\bbL_{\tilde{\mcS}})$ is the non-positive $\tilde{\gamma}$-weight summand of $(\ev_1 \circ \spfilt)^\ast(\bbL_{\mcD})$ by \Cref{L:stratum_cotangent}, but $(\ev_1 \circ \spfilt)^\ast(\bbL_{\mcD})$ is already nonpositively graded by the calculations above, so $\ev_1^\ast (\bbL_{\mcD}) \to \bbL_{\tilde{\mcS}}$ is an isomorphism after pulling back under $\spfilt$. But the only open substack of $\tilde{\mcS}$ containing the image of $\spfilt$ is the whole stack $\tilde{\mcS}$, so $\ev_1^\ast(\bbL_{\mcD}) \to \bbL_{\tilde{\mcS}}$ is an isomorphism.
	\end{itemize}	
	
	We now claim that $\ev_1 : \widetilde{\mcS} \to \mcD$ is an isomorphism. Because we have shown that $\ev_1$ is \'etale, it suffices to show that it is bijective on geometric points. Let $k$ be an algebraically closed field. A $\Spec(k)$-point of $\mcD$ corresponds to a one-dimensional vector space $L$, an element $s : k \to L$, a morphism $x : \Spec(k) \to \mcS$, and if $s = 0$ a morphism $\Spec(k \oplus L^{-1}[1]) \to \mcZ$ over $\mcS$, where $k \oplus L^{-1}[1]$ denotes the free algebra on $L^{-1}$ in homological degree $1$.
	
	Because $\mcS$ is a $\Theta$-stratum with center $\mcZ$, there is a unique morphism $f : \Theta_k \to \mcS$ along with $f(1) \cong x$ such that restriction of $f$ to $\{0\}/\bbG_m$ corresponds to a point of $\mcZ \subset \Grad(\mcS)$. $L \otimes_k \mcO_{\Theta}\langle -1 \rangle$ is the unique invertible sheaf extending $L$ at the point $1 : \Spec(k) \to \Theta_k$ and whose weight at $0 \in \Theta_k$ is $1$, and the section $s$ extends uniquely to a section $t \cdot s$ of $L \otimes_k \mcO_{\Theta}\langle -1 \rangle$. If $s \neq 0$, then $V(s) = \{0\}/\bbG_m$, and as already mentioned there is a unique morphism $\{0\}/\bbG_m \to \mcZ$ over the morphism $f : \Theta_k \to \mcS$. Finally, if $s = 0$, then the morphism $\Spec(k \oplus L^{-1}[1]) \to \mcZ$ corresponds to a lift of $x$ to a $k$-point $\tilde{x} : \Spec(k) \to \mcZ$ and a point of $\Map(\tilde{x}^\ast(\bbL_{\mcZ/\mcS}),L^{-1}[1])$. Now $V(t \cdot s) = V(0) = \Spec_{\Theta_k} (\mcO_\Theta \oplus L^{-1} \otimes_k \mcO_{\Theta} \langle 1 \rangle[1])$, where $\mcO_\Theta \oplus L^{-1} \otimes_k \mcO_{\Theta} \langle 1 \rangle[1]$ is the free algebra on the locally free sheaf $L^{-1} \otimes_k \mcO_\Theta \langle 1 \rangle$ in homological degree $1$. The point $\tilde{x}$ extends uniquely to a morphism $\tilde{f} : \Theta_k \to \mcZ$ lifting $f : \Theta_k \to \mcS$, and we consider the restriction to $\{1\} \in \Theta_k$,
	\begin{equation}\label{E:small_quantization_theorem}
		\Map_{\Theta_k}(\tilde{f}^\ast(\bbL_{\mcZ/\mcS}),L^{-1} \otimes_k \mcO_{\Theta} \langle 1 \rangle[1]) \to \Map_k(\tilde{x}^\ast(\bbL_{\mcZ/\mcS}),L^{-1}).
	\end{equation}
	We regard $\{0 \}/ \bbG_m \hookrightarrow \Theta_k$ as a $\Theta$-stratum with cocharacter $\lambda(t) = t^{-1}$ and semistable locus $(\bbA^1_k \setminus 0) /\bbG_m$, and with respect to this cocharacter $\tilde{f}^\ast(\bbL_{\mcZ/\mcS})|_{\{0\}} \in \QC(B\bbG_m)^{\geq 2}$ and $L^{-1} \otimes_k \mcO_{\Theta} \langle 1 \rangle[1] \in \QC(\Theta_k)^{<2}$. Therefore \eqref{E:small_quantization_theorem} an isomorphism by the quantization commutes with reduction theorem \cite[Prop.~2.1.4]{Halpern-Leistner2021DerivedConjecture}, which completes the proof that our original $k$ point of $\mcD$ lifts uniquely to a $k$-point of $\tilde{\mcS}$.

    The identification $\widetilde{\mcS} \cong \mcD_{\mcZ/\mcS}$ equips $\mcD_{\mcZ/\mcS}$ with a $\Theta$-action coming from the canonical $\Theta$-action on $\Filt(\mcD_{\mcZ/\mcS})$. $\Theta$ also acts on $\mcZ \times \Theta$ via the tautological action on the right factor and the $\Theta$ action on $\mcZ$ via $\gamma$, and the morphism $\tilde{\spfilt} : \mcZ \times \Theta \to \mcD_{\mcZ/\mcS}$ is $\Theta$-equivariant. It follows from \cite[Prop.~1.1.2(4)]{Halpern-Leistner2021DerivedConjecture} that $\tilde{\spfilt}_\ast$ maps $\QC(\mcZ \times \Theta)^{<w}$ to $\QC(\mcD_{\mcZ/\mcS})^{<w}$ and $\tilde{\spfilt}^\ast$ maps $\QC(\mcD_{\mcZ/\mcS})^{<w}$ to $\QC(\mcZ \times \Theta)^{<w}$ for all $w$. Concretely, $\QC(\mcZ \times \Theta)$ corresponds the category of diagrams $\cdots \to F_{n+1} \to F_{n} \to \cdots$ in $\QC(\mcZ)$, and $\QC(\mcZ \times B\bbG_m)^{<w}$ is the full subcategory of diagrams where $F_n^v = 0$ whenever $v + n \geq w$, where $F_n^v$ denotes the weight $v$ summand with respect to the $B\bbG_m$-factor of $\mcZ \times B\bbG_m$. In particular, for $F \in \QC(\mcZ)^{<w}$ with respect to the $\gamma$-action, $F \boxtimes \mcO_{\Theta} \in \QC(\mcZ \times \Theta)^{<w}$ with respect to the $\tilde{\gamma}$-action, and hence $\tilde{F} \in \QC(\mcZ \times \Theta)^{<w}$. This means that for any $v$, the filtration of the weight-$v$ summand
    \[
    \cdots \to F_1^v \to F_0^v \cong f^\ast(f_\ast(F))^v
    \]
    has $F^v_n \cong 0$ for $n \geq w-v$. This implies $\radj{v}(F_n) \cong 0$ for $n \geq w-v$. 
\end{proof}

\begin{ex}
	To illustrate the last proof, suppose $\mcS = \Spec(A) / \bbG_m$ for a non-positively graded algebra and $\mcZ = \Spec(A/I_-) \times B\bbG_m$, where $I_-$ is the ideal generated by homogeneous elements of negative weight. Then $\mcD_{\mcZ / \mcS} \cong \Spec(R(I_-))/\bbG_m^2$, where $R(I_-) \cong \bigoplus_{k<0} A \cdot t^{-k} \oplus \bigoplus_{k \geq 0} I_-^k \cdot t^{-k} \subset A[t^{\pm 1}]$ is the bigraded algebra where $t$ has weight $(0,-1)$ and an element of $A_k$ has weight $(k,0)$. The idea of the previous proof is that we can find a cocharacter of $\bbG_m^2$ with respect to which $R(I_-)$ is non-positively graded, and the ideal generated by negative weight elements is $\bigoplus_{k<0} A \cdot t^{-k} \oplus I_- \cdot t^0 \oplus \bigoplus_{k>0} I_-^k \cdot t^{-k}$, so that the fixed locus is again $\Spec(A/I_-)$. This is achieved by the cocharacter $\tilde{\gamma}(t) = (t^2 , t)$.
\end{ex}

\subsection{Baric completion in general}

We say that a baric structure is \emph{right-complete} if $\forall F \in \mcC$, the morphisms $\radj{w}(F) \to F$ realize $F$ as $\colim_{w \to -\infty} \radj{w}(F)$.

\begin{defn}[Baric completion]\label{D:baric_completion}
	Given a stable subcategory $\mcA \subset \mcC$, we define the baric completion
	\[
		\mcA^\wedge_{\beta} = \left\{K \in \mcC \left| \begin{array}{c} \radj{w}(K) \in \mcA \text{ for all }w \in \bbZ \\ K \in \mcC^{<w} \text{ for some } w \in \bbZ \end{array} \right.\right\}.
	\]
	We also define $\mcC^{\loc} := \bigcap_w \mcC^{<w}$, $\mcA^{\loc} := \mcC^{\loc} \cap \mcA$, and
	\[
		\mcA^{\nil} := \{ F \in \mcA^\wedge_\beta | F \cong \underset{w \to -\infty}{\colim} \radj{w}(F) \}.
	\]
	The superscripts stand for ``local" and ``nilpotent."
\end{defn}

\begin{ex} \label{ex:baric_structure_representations}
	Let $G$ be a linearly reductive group over a field $k$. Because $\QC(BG)$ is semisimple, any function from the set of irreducible representations of $G$ to $\bbZ$ defines a baric decomposition, where $\QC(\mcX)^{< w}$ is the subcategory generated by representations assigned value $< w$. There is an interesting special case, however, that arises in examples: Let $M$ be weight lattice of a maximal torus $T_{\bar{k}} \subset G_{\bar{k}}$, where $\bar{k}$ is an algebraic closure of $k$. Let $\sigma \subset M_\bbQ$ be a Weyl group invariant rational polyhedral cone of full dimension and with interior $\sigma^\circ$, and let $v \in - \sigma^\circ$. We define
	\[
		\QC(BG)^{<w} = \left\{ V \in \QC(BG) \mid \forall \chi \in M,  H^\ast(V_{\bar{k}})_\chi = 0 \text{ unless } \chi \in w v + \sigma^\circ \right\},
	\]
    where $(-)_\chi$ denotes the weight $\chi$ summand of a representation of $T_{\bar{k}}$. Then $\Perf(BG)^\wedge_{\beta} \subset \QC(BG)$ is the full subcategory of complexes of representations $V$ such that after base change to $G_{\bar{k}}$, there is a $u \in M$ such that all $\chi \in M$ with weight space $V_\chi \neq 0$ lie in $u + \sigma$, and for any $u \in M$, the subspace $\bigoplus_{\chi \in u - \sigma} V_\chi \subset V_{\bar{k}}$ is finite dimensional. In this case $\Perf(BG)^\wedge_\beta$ is closed under tensor product.
\end{ex}

\begin{lem}
	For any $F \in \mcA^\wedge_\beta$, $\underset{w \to -\infty}{\colim} \radj{w}(F) \in \mcC$ lies in $\mcA^\wedge_\beta$.
\end{lem}
\begin{proof}
	For any $v$, $\radj{v}(\colim_{w \to -\infty} \radj{w}(F)) = \radj{v}(F) \in \mcA$. Also, if $u \in \bbZ$ is such that $\radj{u}(F) \cong 0$, then \[\radj{u}(\underset{w \to -\infty}{\colim}(\radj{w}(F))) \cong \underset{w \to -\infty}{\colim} \radj{w}(\radj{u}(F)) \cong 0.\]
\end{proof}

\begin{lem}\label{L:decompositions_of_baric_completion}
	The baric structure on $\mcC$ restricts to $\mcA^\wedge_{\beta}$, $\mcA^\wedge_{\beta} \cap \mcA$, and $\mcA^{\nil}$. The resulting baric structure on $\mcA^{\nil}$ is right-complete. We have semiorthogonal decompositions
	\[
		\mcA^\wedge_{\beta} = \langle \mcC^{\loc}, \mcA^{\nil} \rangle \quad \text{and} \quad \mcA^\wedge_{\beta} \cap \mcA = \langle \mcA^{\loc}, \mcA^{\nil} \rangle.
	\]
\end{lem}

\begin{proof}
	The definitions imply that $\cofib(\colim_{w \to -\infty} \radj{w}(F) \to F) \in \mcC^{\loc}$, and $R\Hom(\mcA^{\nil},\mcC^{\loc}) \cong 0$.
\end{proof}

\begin{lem}\label{L:K_theory_graded}
	Let $\mcC$ be a stable $\infty$-category, and consider a diagram $\cdots \to F_{i+1} \to F_{i} \to \cdots$ such that $F_i \simeq 0$ for $i \gg 0$. Let $G_i := \cofib(F_{i+1} \to F_i)$ for all $i$. If $F := \colim_i F_i$ exists in $\mcC$ and $\bigoplus_{i \leq a} \cofib(F_i \to F)$ exists in $\mcC$ for some $a\in \bbZ$, then $\bigoplus_i G_i$ exists, and $[F] = [\bigoplus_i G_i]$ in $K_0(\mcC)$.
\end{lem}
\begin{proof}
	If $\bigoplus_{i \leq a} \cofib(F_i \to F)$ exists for some $a$, then it exists for any $a$. So we may choose $a \gg 0$ so that $F_i \cong 0, \forall i \geq a$. The claim of the lemma follows from the exact triangle
	\[
		F \oplus \bigoplus_{i \leq a} \cofib(F_i \to F) \to \bigoplus_{i \leq a} \cofib(F_i \to F) \to \bigoplus_{i < a} G_i[1] = \bigoplus_{i \in \bbZ} G_i[1],
	\]
	where the first arrow is the direct sum of the canonical morphisms $\cofib(F_i \to F) \to \cofib(F_{i-1} \to F)$ for $i \leq a$ and the isomorphism $F \to \cofib(F_a \to F)$.
\end{proof}

\begin{lem} \label{L:baric_splitting_K_theory}
	Baric projection defines a functor
	\[
		\prod_w \beta^w \colon \mcA^{\nil} \to \bigoplus_{w \geq 0} \mcA^w \oplus \prod_{w<0} \mcA^w
	\]
	that induces an isomorphism on $K_0(-)$.
\end{lem}
\begin{proof}
	This amounts to the claim that the functor $E \mapsto \bigoplus_{w} \beta^w(E)$ induces the identity homomorphism on $K_0(\mcA^\nil)$. We know that $\colim_{w \to -\infty} \radj{w}(E) \cong E$, and the $\beta^w(F) = \cofib(\radj{w+1}(F) \to \radj{w}(F))$. Note that $\bigoplus_{i\leq 0} \ladj{i}(F)$ exists in $\mcA^\nil$, so this follows from \Cref{L:K_theory_graded}.
\end{proof}

\begin{ex}
    Let $G$ be a reductive group over a field $k$, and suppose $T \subset G$ is a split maximal torus with weight lattice $M$. Then $K_0(\Perf(BT)) \cong \bbZ[M]$ and $K_0(\Perf(BG)) \cong \bbZ[M]^W \subset \bbZ[M]$, where $W \subset G$ is the Weyl group of $T$.
    
    If one chooses a pair $(\sigma,v)$ as in \Cref{ex:baric_structure_representations}, then $\bbZ[\sigma \cap M] \subset \bbZ[M]$ is the coordinate ring of an affine toric variety, and one has an isomorphism as rings
    \[
    \begin{array}{rl} K_0 \left(\Perf(BT)_\beta^\wedge \right) &\cong \widehat{\bbZ[\sigma \cap M]} \otimes_{\bbZ[\sigma \cap M]} \bbZ[M], \text{ and }\\
    K_0\left(\Perf(BG)^\wedge_\beta \right) &\cong \left(\widehat{\bbZ[\sigma \cap M]} \otimes_{\bbZ[\sigma \cap M]} \bbZ[M]\right)^W,
    \end{array}
    \]where the completion is with respect to the ideal in $\bbZ[\sigma \cap M]$ generated by $\sigma^\circ \cap M$. These calculations follow from \Cref{L:baric_splitting_K_theory}, which gives a description of the $K$-theory of the baric completions as subgroups of the formal product group $\prod_{\chi \in M} \bbZ \cdot e^\chi$
\end{ex}

\subsection{Baric completion in our context} \label{S:baric_completion_perf}

We now consider a derived stack $\mcS$ with a weak $\Theta$-action, which equips $\QC(\mcS)$ with a baric structure as above. We will consider the baric completion $\Perf(\mcS)^{\wedge}_\beta$. It follows from \cite[Prop.~1.1.2(2)]{Halpern-Leistner2021DerivedConjecture} that $(\Perf(\mcS)^{\wedge}_\beta)^{\nil}=\Perf(\mcS)^{\wedge}_\beta$, i.e., the baric structure on $\Perf(\mcS)^{\wedge}_\beta$ is right-complete.

\begin{lem}\label{L:perf_completion_functoriality}
	$E \in \QC(\mcS)$ lies in $\Perf(\mcS)^\wedge_\beta$ if and only if $E \in \QC(\mcS)^{<\infty}$ and $\spfilt^\ast(E) \in \Perf(\mcZ)^\wedge_\beta$, and $F \in \QC(\mcZ)$ lies in $\Perf(\mcZ)^\wedge_\beta$ if and only if $\gr^\ast(F) \in \Perf(\mcS)^\wedge_\beta$.
\end{lem}
\begin{proof}
	The functors $\spfilt^\ast$ and $\gr^\ast$ preserve the subcategories $\Perf(-)^\wedge_\beta$ because they commute with baric truncation \cite[Prop.~1.1.2(4)]{Halpern-Leistner2021DerivedConjecture} and preserve $\Perf(-)$. This immediately implies that if $F \in \QC(\mcZ)$ and $\gr^\ast(F) \in \Perf(\mcS)^\wedge_\beta$, then $F \cong \spfilt^\ast(\gr^\ast(F)) \in \Perf(\mcZ)^\wedge_\beta$.

	On the other hand, for $E \in \QC(\mcS)^{<\infty}$, the canonical isomorphism $E \cong \colim_{w \to -\infty} \radj{w}(E)$ gives a bounded below filtration of $E$ whose graded pieces are $\beta^w(E) \cong \gr^\ast(\beta^w(\spfilt^\ast(E)))$ by \Cref{L:constant_weight_category}, which implies that if $\spfilt^\ast(E) \in \Perf(\mcZ)^\wedge_\beta$, then $E \in \Perf(\mcS)^\wedge_\beta$.
\end{proof}

\begin{lem} \label{L:K_theory_perf_completion}
	The functors $\spfilt^\ast$ and $\gr^\ast$ provide inverse equivalences $K_0(\Perf(\mcS)^\wedge_\beta) \cong K_0(\Perf(\mcZ)^\wedge_\beta)$.
\end{lem}
\begin{proof}
	This follows from \Cref{L:baric_splitting_K_theory} and \Cref{L:constant_weight_category}.
\end{proof}

\begin{ex}\label{EX:perf_completion_Z}
	Every object in $\QC(\mcZ)$ splits as a direct sum $E \cong \bigoplus \beta^w(E)$. Using this, one can show that $\Perf^\wedge_\beta(\mcZ) \subset \QC(\mcZ)$ is the full subcategory of $E$ such that $\beta^w(E) \in \Perf(\mcZ), \forall w$, and $\beta^w(E) \cong 0$ for all $w\gg 0$. In addition, because $K_0(-)$ commutes with infinite products \cite{Kasprowski2019AlgebraicComplexes},
	\[K_0(\Perf(\mcZ)^{\wedge}_\beta) \cong \bigoplus_{w>0} K_0(\Perf(\mcZ)^w) \oplus \prod_{w \leq 0} K_0(\Perf(\mcZ)^w).\]
    Now suppose $L \in \QC(\mcZ)^w$ is an invertible sheaf for some $w<0$. Tensor product with $L$ gives an isomorphism $\Cycles(\mcZ/\mcB)^a \cong \Cycles(\mcZ/\mcB)^{a+w}$ for all $a$. Also, if we let $\mcZ_{\rm{rig}} := \mcZ \times_{B\bbG_m} \rm{pt}$ be the fiber of the morphism $\mcZ \to B\bbG_m$ classified by $L$, then pullback gives an equivalence $\QC(\mcZ_{\rm{rig}}) \cong \bigoplus_{a=0}^{w-1} \QC(\mcZ)^a$. Combining all of this with \Cref{L:K_theory_perf_completion} gives an isomorphism
	\[
	K_0(\Perf(\mcS)^\wedge_\beta) \cong K_0(\Perf(\mcZ_{\rm{rig}}))(\!(u)\!),
	\]
	where the action of $u$ can be identified with $\gr^\ast(L) \otimes (-)$ on the left-hand-side.
\end{ex}

\begin{lem}\label{L:perf_completion_symmetric_monoidal}
	For $E \in \Perf(\mcS)^\wedge_\beta$ and $F \in \QC(\mcS)$, $\hwt(E \otimes F) \leq \hwt(E) + \hwt(F)$, and $E \otimes (-)$ preserves $\Perf(\mcS)^\wedge_\beta$. In particular, $\Perf(\mcS)^\wedge_\beta \subset \QC(\mcS)$ is a symmetric monoidal subcategory containing $\Perf(\mcS)$. The same holds for $\mcZ$.
	 \end{lem}
\begin{proof}
	Suppose  $E \in \QC(\mcX)^{<u}$ and $F \in \QC(\mcX)^{<v}$, with $E \in \Perf(\mcS)^\wedge_\beta$. Using \cite[Prop.~1.1.2(6)]{Halpern-Leistner2021DerivedConjecture} one can show that $\Perf(\mcS)^{<a} \otimes \QC(\mcS)^{<b} \subset \QC(\mcS)^{<a+b-1}$ for any $a,b\in \bbZ$. Because $\QC(\mcS)^{<u+v-1}$ is closed under filtered colimits, this implies that $E \otimes F \cong \colim_{w \to -\infty} \radj{w}(E) \otimes F \in \QC(\mcS)^{<u+v-1}$. If in addition $F \in \Perf(\mcS)^\wedge_\beta$, then the same identity implies that $\forall w \in \bbZ$,
	\[\radj{w}(E \otimes F) \cong \radj{w}(\radj{w-v+1}(E) \otimes \radj{w-u+1}(F)) \in \Perf(\mcS),\]
	hence $E \otimes F \in \Perf(\mcS)^\wedge_\beta$. The same arguments also apply to $\Perf(\mcZ)^\wedge_\beta$, or one can deduce the claim from the explicit description in \Cref{EX:perf_completion_Z}.
\end{proof}

\begin{lem}\label{L:fundamental_class_perf_complete}
	If $\bbL_{\mcS/\mcR} \in \Perf(\mcS)$, then the complex $\spfilt_\ast(\mcO_\mcZ) \in \QC(\mcS)$ lies in $\Perf(\mcS)^\wedge_\beta$, and $\gr_\ast(\mcO_\mcS) \in \Perf(\mcZ)^\wedge_\beta$.
\end{lem}
\begin{proof}
	It follows from adjunction that $\spfilt_\ast(\mcO_\mcZ) \in \QC(\mcS)^{<\infty}$, so it suffices by \Cref{L:perf_completion_functoriality} to show that $\spfilt^\ast(\spfilt_\ast(\mcO_\mcZ)) \in \Perf(\mcZ)^\wedge_\beta$. $\bbL_{\mcZ/\mcS} \cong \ladj{0}(\spfilt^\ast(\bbL_{\mcX/\mcR}))[1]$ by \Cref{L:stratum_cotangent}, so $\Sym^n(\bbL_{\mcZ/\mcS}) \in \APerf(\mcZ)^{<n}$ for all $n \geq 0$. By \Cref{T:push_pull_filtration_stratum}, for any $w \in \bbZ$, $\radj{w}(\spfilt^\ast(\spfilt_\ast(\mcO_\mcZ)))$ has a finite filtration whose associated graded pieces are $\radj{w}(\Sym^n(\bbL_{\mcZ/\mcS}))$. If $\bbL_{\mcS/\mcR} \cong \ladj{1}(\ev_1^\ast(\bbL_{\mcX/\mcR}))$ is perfect, then so is $\Sym^n(\bbL_{\mcZ/\mcS})$ for all $n\geq 0$. Therefore $\radj{w}(\spfilt^\ast(\spfilt_\ast(\mcO_\mcZ))) \in \Perf(\mcZ)$ for all $w$, which proves that $\spfilt_\ast(\mcO_\mcZ) \in \Perf(\mcS)^\wedge_\beta$ by \Cref{L:perf_completion_functoriality}.

	The proof that $\gr_\ast(\mcO_\mcS) \in \Perf(\mcZ)^\wedge_\beta$ is identical, using the fact that $\spfilt^\ast \bbL_{\mcS/\mcZ} \cong \bbL_{\mcZ/\mcS}[-1]$ and hence $\bbL_{\mcS/\mcZ} \in \QC(\mcS)^{<0}$.
\end{proof}

\begin{lem}\label{L:perf_completion_functoriality_2}
	If $\bbL_{\mcS/\mcR} \in \Perf(\mcS)$, then the functors $\spfilt_\ast : \QC(\mcZ) \to \QC(\mcS)$ and $\gr_\ast : \QC(\mcS) \to \QC(\mcZ)$ preserve the subcategories $\Perf(-)^\wedge_\beta$.
\end{lem}
\begin{proof}
	The projection formula implies that
	\[
		\spfilt_\ast(E) \cong \gr^\ast(E) \otimes \spfilt_\ast(\mcO_\mcZ).
	\]
	The claim for $\spfilt_\ast$ then follows from \Cref{L:perf_completion_functoriality}, \Cref{L:fundamental_class_perf_complete}, and \Cref{L:perf_completion_symmetric_monoidal}.

	On the other hand, $\gr_\ast( \QC(\mcS)^{<w}) \subset \QC(\mcZ)^{<w}$ by adjunction, because $\gr^\ast(\QC(\mcZ)^{\geq w}) \subset \QC(\mcS)^{\geq w}$. Therefore, for $E \in \Perf(\mcS)^\wedge_\beta$,
	\[
		\radj{w}(\cofib(\gr_\ast(\radj{w}(E)) \to \gr_\ast(E))) \cong 0,
	\]
	so to prove that $\gr_\ast(E) \in \Perf(\mcZ)^\wedge_\beta$, it suffices to prove this for $\gr_\ast(\radj{w}(E))$ for every $w$. By \Cref{L:constant_weight_category}, $\radj{w}(E)$ is a finite sequence of extensions of objects of the form $\gr^\ast(F)$ for $F\in \Perf(\mcZ)$, so it suffices to show that $\gr_\ast(\gr^\ast(F)) \cong F \otimes \gr_\ast(\mcO_\mcS) \in \Perf(\mcZ)^\wedge_\beta$, which follows from \Cref{L:fundamental_class_perf_complete} and \Cref{L:perf_completion_symmetric_monoidal}.
\end{proof}

\subsection{The Euler class of a complex}\label{S:euler_class}

Given a vector bundle $E$ over a space $X$, the K-theoretic Euler class is defined as $e(E) := [\Sym(E^\ast[1])] \in K_0(X)$. Using the perfect pairing on the exterior algebra, one can also write this as $e(E) = [\Sym(E[1])\otimes \det(E[1]) [\rank(E)]]$. For short exact sequences of vector bundles $0 \to E \to F \to G \to 0$ one has $e(F) = e(E) \cdot e(G)$, but $e(-)$ does not extend to a function on $K$-theory classes.

Now let $\mcS$ be a stack with a weak $\Theta$-action as in the previous section. Given $E \in \Perf(\mcS)$ such that $\beta^0(E) \cong 0$, we can define the \emph{Euler class}
\[
	\begin{array}{rll}
		e(E) & := \Sym(E^\ast[1])                          & \text{if } E \in \Perf(\mcS)^{\geq 1} \\
		     & := \Sym(E[1]) \otimes \det(E[1])[-\rank(E)] & \text{if }E \in \Perf(\mcS)^{<0}      \\
		     & := e(\ladj{0}(E)) \cdot e(\radj{1}(E))      & \text{ in general}
	\end{array}
\]
Note that $E \in \Perf(\mcS)^{\geq 1} \Rightarrow E^\ast \in \Perf(\mcS)^{<0}$, and $\Sym^n(\Perf(\mcS)^{<0}) \subset \Perf(\mcS)^{<n-1}$ by \cite[Lem.~1.5.4]{Halpern-Leistner2021DerivedConjecture}. Therefore the complex $e(E) \in \QC(\mcS)$ actually lies in $\Perf(\mcS)^{\wedge}_\beta$. Because $\Perf(\mcS)_{\beta}^\wedge$ is a symmetric monoidal category, by \Cref{L:perf_completion_symmetric_monoidal}, its $K$-theory is naturally a ring. In fact, $e(E)$ is always a unit in this ring.

\begin{lem}\label{L:euler_class}
	Let $\Perf(\mcS)^{\rm{mov}} := \{E \in \Perf(\mcS) | \beta^0(E) \cong 0\}$. Then $e(-)$ defines a group homomorphism
	\[K_0(\Perf(\mcS)^{\rm{mov}}) \to K_0(\Perf(\mcS)_{\beta}^\wedge)^\times.\]
\end{lem}
\begin{proof}
    We must show that $e(F) = e(E) \cdot e(G)$ in $K_0(\Perf(\mcS)^\wedge_\beta)$ for any exact triangle $E \to F \to G$ in $\Perf(\mcS)^{\rm{mov}}$. The two-step filtration of $F$ corresponding to this exact triangle induces a finite filtration on $\Sym^n(F)$ for any $n \geq 0$ whose associated graded is $\Sym^n(E \oplus G)$, so $[\Sym^n(F)] = [\Sym^n(E \oplus G)] \in K_0(\Perf(\mcS))$. If $E,F,$ and $G$ are in $\Perf(\mcS)^{<0}$, then $\bigoplus_{n \geq a} \Sym^n(F)$ and $\bigoplus_{n \geq a} \Sym^n(E \oplus G)$ lie in $\Perf(\mcS)^\wedge_\beta$, so \Cref{L:K_theory_graded} implies that $[\Sym(F)] = [\Sym(E \oplus G)] = [\Sym(E) \otimes \Sym(G)]$ in $\Perf(\mcS)^\wedge_\beta$.
\end{proof}

Using this lemma, the expression $e(\bbN_{\mcZ_\alpha/\mcX})^{-1}$ in the non-abelian localization formula will be interpreted as the inverse in $K_0(\Perf(\mcZ_\alpha)^\wedge_\beta)$.

\section{The non-abelian localization theorem}

\subsection{Highest weight cycles}

Now we let $\rho : \mcX \to \mcB$ be as in \Cref{H:lfp}, and we also assume that $\mcX$ is quasi-compact. We will consider a single closed $\Theta$-stratum $\mcS \subset \Filt(\mcX)$ relative to $\mcB$. (See \Cref{D:relative_theta_stratification}.)

\begin{defn}[Highest weight cycles]\label{D:highest_weight_cycles}
	Given the baric structure on $\QC(\mcX)$ above, we define the subcategory of \emph{highest weight cycles} on $\mcZ$, $\mcS$, and $\mcX$:
	\begin{align*}
        \Cycles(\mcZ/\mcB)^{<\infty} &:= \Cycles(\mcZ/\mcB) \cap \QC(\mcZ)^{<\infty} \\
        \Cycles(\mcS/\mcB)^{<\infty} &:= \Cycles(\mcS/\mcB) \cap \QC(\mcS)^{<\infty} \\
        \Cycles(\mcX/\mcB)^{<\infty} &:= \left\{ F \in \QC(\mcX)^{<\infty} \mid F|_{\mcX \setminus \mcS} \in \Cycles(\mcX\setminus \mcS) \text{ and }\ev_1^!(F) \in \Cycles(\mcS/\mcB)\right\}
	\end{align*}
    We also let $\Cycles(-)^{<v} := \Cycles(-)^{<\infty} \cap \QC(-)^{<v}$ for any $v \in \bbZ$.
\end{defn}

Note that despite the notation, it is not obvious that $\Cycles(\mcX/\mcB)^{<\infty} \subset \Cycles(\mcX/\mcB)$, but we will show this in \Cref{L:highest_weight_cycles_are_cycles}.

\begin{ex}\label{EX:highest_weight_cycles_on_Z}
	The category $\Cycles(\mcZ/\mcB)^{<\infty}$ is simpler than the others. Under the isomorphism $\QC(\mcZ) \cong \bigoplus_{w \in \bbZ} \QC(\mcZ)^w$, a complex $E = \bigoplus_w E^w$ lies in $\Cycles(\mcZ/\mcB)^{<\infty}$ if and only if $E^w \cong 0$ for $w \gg 0$ and each $E^w \in \Cycles(\mcZ/\mcB)$. In fact, the second condition can be replaced with the condition that $(\rho|_\mcZ)_\ast(P \otimes E^w) \in \DCoh(\mcB)$ for any $w \in \bbZ$ and $P \in \Perf(\mcZ)^{-w}$. The same considerations as in \Cref{EX:perf_completion_Z} show that if $L \in \QC(\mcZ)^w$ for some $w<0$ is an invertible sheaf and $\mcZ_{\rm{rig}}:= \mcZ \times_{B\bbG_m} \rm{pt} \cong \rm{Tot}(L) \setminus 0$ is the corresponding rigidification, then
	\[K_0(\Cycles(\mcZ/\mcB)^{<\infty}) \cong K_0(\Cycles(\mcZ_{\rm{rig}}))(\!(u)\!),\]
	where $u = [L] \in K_0(\Perf(\mcZ))$.
\end{ex}

\begin{lem}\label{L:vanishing_rho}
	For $F \in \Perf(\mcS)^{<w}$, $(\rho|_\mcS)_\ast(F \otimes (-)) \cong (\rho|_\mcS)_\ast(F \otimes \radj{1-w}(-))$ as functors $\QC(\mcS) \to \QC(\mcB)$.
\end{lem}
\begin{proof}
	\Cref{L:perf_completion_symmetric_monoidal} implies that $F \otimes \QC(\mcS)^{<1-w} \subset \QC(\mcS)^{<0}$. We claim that $(\rho|_{\mcS})_\ast$ vanishes on $\QC(\mcS)^{<0}$. By base change along a smooth cover $\Spec(A) \to \mcB$, it suffices to prove this when $\mcB$ is affine -- this is where we use that $\mcS$ is a $\Theta$-stratum \emph{relative to $\mcB$}. When $\mcB$ is affine, the vanishing follows from semiorthogonality between $\mcO_\mcS \in \QC(\mcS)^{\geq 0}$ and $\QC(\mcS)^{<0}$. Finally, the vanishing implies the lemma, because $(\rho|_\mcS)_\ast(F \otimes \ladj{1-w}(-)) \cong 0$.

\end{proof}

\begin{lem}\label{L:highest_weight_cycles_module_category}
	 For any $E \in \Perf(-)^\wedge_\beta$, $E \otimes(-)$ preserves $\Cycles(\mcS/\mcB)^{<\infty}$. The analogous claim holds for $\mcZ$.
\end{lem}
\begin{proof}
	For any $P \in \Perf(\mcS)$ and $F \in \Cycles(\mcS/\mcB)^{<\infty}$, $P \otimes F \in \QC(\mcS)^{<v}$ for some $v \in \bbZ$, so \Cref{L:vanishing_rho} implies that $(\rho|_\mcS)_\ast(P \otimes E \otimes F) \cong (\rho|_\mcS)_\ast(P \otimes \radj{1-v}(E) \otimes F)$, which lies in $\DCoh(\mcB)$ because $P \otimes \radj{w}(E) \in \Perf(\mcS)$ by hypothesis. \Cref{L:perf_completion_symmetric_monoidal} also implies that $E \otimes F \in \QC(\mcS)^{<\infty}$, so $E \otimes F \in \Cycles(\mcS/\mcB)^{<\infty}$. The same proof shows that $\Perf(\mcZ)^\wedge_\beta \otimes \Cycles(\mcZ/\mcB)^{<\infty} \subset \Cycles(\mcZ/\mcB)^{<\infty}$.
\end{proof}

For notational convenience, we will also denote $\QC(\mcX \setminus \mcS)^{<\infty} := \QC(\mcX \setminus \mcS)$ and $\Cycles(\mcX \setminus \mcS)^{<\infty} := \Cycles(\mcX \setminus \mcS)$.

\begin{prop} \label{P:functoriality_highest_weight_cycles}
	All of the functors in the following diagram preserve the subcategories $\Cycles(-)^{<\infty} \subseteq \QC(-)$:
	\[
		\xymatrix{\QC(\mcZ) \ar@/^1pc/@<1ex>[r]^-{\spfilt_\ast} \ar@/_1pc/[r]^-{\gr^\ast} & \QC(\mcS) \ar@/^1pc/@<1ex>[l]^{\gr_\ast} \ar@/_1pc/[l]^-{\spfilt^\ast} \ar@<.5ex>[r]^-{\ev_{1\ast}} & \ar@<.5ex>[l]^-{\ev_1^!} \QC(\mcX) \ar@<-.5ex>[r]_-{i^\ast} & \ar@<-.5ex>[l]_-{i_\ast} \QC(\mcX \setminus \mcS) }
	\]
\end{prop}
\begin{proof}
	The claims for $\ev_1^!$, $i^\ast$, and $i_\ast$ are just restating the definitions. The functors $\spfilt_\ast$ and $\gr_\ast$ preserve $\Cycles(-)$ automatically, and they preserve $\QC(-)^{<\infty}$ because their left adjoints commute with baric truncation, and thus preserve $\QC(-)^{\geq w}$.

	To show that $\ev_{1\ast}$ preserves $\Cycles(-)^{<\infty}$, it suffices to show that $\ev_1^! \circ \ev_{1\ast}$ maps $\Cycles(\mcS/\mcB)^{<\infty}$ to $\Cycles(\mcS/\mcB)$, because $\ev_{1\ast}$ and $\ev_1^!$ both preserve $\QC(-)^{<\infty}$. For any $E \in \Cycles(\mcS/\mcB)^{<w}$, \Cref{L:vanishing_rho} implies that $(\rho|_{\mcS})_\ast(F \otimes \ev_1^! \circ \ev_{1\ast}(E)) \cong (\rho|_{\mcS})_\ast(F \otimes \radj{1-w}(\ev_1^! \circ \ev_{1\ast}(E)))$, and by \Cref{P:push_pull_filtration_closed_immersion} the latter has a finite filtration whose associated graded pieces are $(\rho|_{\mcS})_\ast(F \otimes \radj{1-w}(\Sym^n(\bbL_{\mcS/\mcX}) \otimes E)) \cong (\rho|_{\mcS})_\ast(F \otimes \Sym^n(\bbL_{\mcS/\mcX}) \otimes E) \in \DCoh(\mcB)$, and hence $(\rho|_{\mcS})_\ast(F \otimes \ev_1^! \circ \ev_{1\ast}(E)) \in \DCoh(\mcB)$ for all $F \in \Perf(\mcS)$.

	The functor $\gr^\ast$ commutes with baric truncation, so it preserves $\QC(-)^{<\infty}$. Using the baric structure on $\mcS$, any $P \in \Perf(\mcS)$ has a finite filtration whose associated graded lies in the essential image of $\gr^\ast$. To show that $\gr^\ast$ preserves $\Cycles(-)$, it therefore suffices to show that $(\rho|_\mcS)_\ast(\gr^\ast(F)) \cong (\rho|_\mcZ)_\ast(F \otimes \gr_\ast(\mcO_\mcS))$ is finite dimensional for any $F \in \Cycles(\mcZ/\mcB)$. This follows from \Cref{L:highest_weight_cycles_module_category} and \Cref{L:fundamental_class_perf_complete}.

	Finally, $\spfilt^\ast$ commutes with baric truncation, and thus preserves $\QC(-)^{<\infty}$. \Cref{L:baric_truncation_preserves_highest_weight_cycles} below -- which depends on the fact established above that $\gr^\ast$ and $\ev_{1\ast}$ preserve $\Cycles(-)^{<\infty}$ but does not depend on the claim for $\spfilt^\ast$ -- reduces the last claim to showing that $\spfilt^\ast(E) \in \Cycles(\mcZ/\mcB)$ for any $E \in \Cycles(\mcS/\mcB)^w$. \Cref{L:constant_weight_category} implies that for any such $E$ and any $P \in \Perf(\mcZ)^w$, we have $R\Hom_\mcZ(P,\spfilt^\ast(E)) \cong R\Hom_\mcS(\gr^\ast(P),E)$, which is finite dimensional. Therefore, $\spfilt^\ast$ preserves $\Cycles(-)^{<\infty}$.
\end{proof}

\begin{prop}\label{L:baric_truncation_preserves_highest_weight_cycles}
	If $K \in \QC(\mcX)^{<\infty}$ and $K|_{\mcX\setminus \mcS} \in \Cycles(\mcX \setminus \mcS)$, then $K \in \Cycles(\mcX/\mcB)^{<\infty}$ if and only if $\beta^w(K) \in \Cycles(\mcX/\mcB)^{<\infty}$ for all $w$. In particular the baric truncation functors $\radj{w}$ and $\ladj{w}$ preserve $\Cycles(\mcX/\mcB)^{<\infty}$.
\end{prop}
\begin{rem}
	The statement of \Cref{L:baric_truncation_preserves_highest_weight_cycles} implies the analogous claim that $K \in \QC(\mcS)^{<\infty}$ lies in $\Cycles(\mcS/\mcB)^{<\infty}$ if and only if $\forall w, \radj{w}(K) \in \Cycles(\mcS/\mcB)$, and therefore the baric structure on $\QC(\mcS)$ preserves $\Cycles(\mcS/\mcB)^{<\infty}$. The same holds for $\Cycles(\mcZ/\mcB)^{<\infty}$. Indeed, these are special cases of \Cref{L:baric_truncation_preserves_highest_weight_cycles} applied to $\mcS$ or $\mcZ$ regarded canonically as a $\Theta$-stratum within itself, so that $\mcS = \mcX$ or $\mcZ=\mcS=\mcX$.
\end{rem}

\begin{proof}[Proof of \Cref{L:baric_truncation_preserves_highest_weight_cycles}]
	We will prove that for $K \in \QC(\mcX)^{<\infty}$, $\ev_1^!(K)
		\in \Cycles(\mcS/\mcB)$ if and only if $\ev_1^!(\beta^w(K)) \in \Cycles(\mcS/\mcB)$ for all $w \in \bbZ$.

	Suppose that $\forall w, \ev_1^!(\beta^{w}(K)) \in \Cycles(\mcS/\mcB)$. This implies $\forall w, \ev_1^!(\radj{w}(K)) \in \Cycles(\mcS/\mcB)$, because $K \in \QC(\mcX)^{<\infty}$. For any $P \in \Perf(\mcS)$, \Cref{L:vanishing_rho} implies that for some $w \ll 0$, $(\rho|_\mcS)_\ast(P \otimes \ev_1^!(K)) \cong (\rho|_\mcS)_\ast(P \otimes \ev_1^!(\radj{w}(K))) \in \DCoh(\mcB)$. Therefore $\ev_1^!(K) \in \Cycles(\mcS/\mcB)$.

	For the converse, it suffices by an inductive argument to show that if $K \in \QC(\mcX)^{<w+1}$ and $\ev_1^!(K) \in \Cycles(\mcS/\mcB)$, then $\ev_1^!(\radj{w}(K)) \in \Cycles(\mcS/\mcB)$. By \Cref{L:constant_weight_category}, $\radj{w}(K) \cong \ev_{1\ast} \circ \gr^\ast(G)$ for some $G \in \QC(\mcZ)^w$. For any $P \in \Perf(\mcZ)^w$,
	\begin{align*}
		R\Hom_\mcZ(P,G) & \cong R\Hom_\mcX(\ev_{1\ast} \circ \gr^\ast(P), \radj{w}(K)) \\
		               & \cong R\Hom_\mcX(\ev_{1\ast} \circ \gr^\ast(P),K)            \\
		               & \cong R\Hom_\mcS(\gr^\ast(P),\ev_1^!(K)),
	\end{align*}
	which is finite dimensional because $\ev_1^!(K) \in \Cycles(\mcS/\mcB)$ by hypothesis. Therefore $G \in \Cycles(\mcZ/\mcB)^w$, so $\radj{w}(K) \cong \ev_{1\ast} \circ \gr^\ast(G) \in \Cycles(\mcX/\mcB)^{<\infty}$ by \Cref{P:functoriality_highest_weight_cycles} and therefore $\ev_1^!(\radj{w}(K)) \in \Cycles(\mcS/\mcB)$ by definition.

	This shows that the baric truncation functors preserve $\Cycles(\mcX/\mcB)^{<\infty}$, because $\radj{w}(K)|_{\mcX \setminus \mcS} \cong 0$ and $K|_{\mcX \setminus \mcS} \cong \ladj{w}(K)|_{\mcX\setminus \mcS}$.
\end{proof}

\begin{lem}\label{L:constant_weight_cycles}
The functors $\gr^\ast : \Cycles(\mcZ/\mcB)^w \to \Cycles(\mcS/\mcB)^w$ and $(\ev_1)_\ast : \Cycles(\mcS/\mcB)^w \to \Cycles(\mcX/\mcB)^w$ are equivalences of categories.
\end{lem}
\begin{proof}
By \Cref{L:constant_weight_category}, these functors are equivalences on $\QC(-)^w$. The inverse of $\gr^\ast$ is $\spfilt^\ast$, and the inverse of $(\ev_1)_\ast$ is $\radj{w}(\ev_1^!(-))$. All four functors preserve $\Cycles(-)^w$ by \Cref{P:functoriality_highest_weight_cycles} and \Cref{L:baric_truncation_preserves_highest_weight_cycles}.
\end{proof}

\begin{cor}\label{C:completion_tensor_property}
    If $E \in \QC(\mcX)$ is such that $\ev_1^\ast(E) \in \Perf(\mcS)^\wedge_\beta$ and $E|_{\mcX \setminus \mcS} \in \Perf(\mcX \setminus \mcS)$, then $\hwt(E \otimes F) \leq \hwt(\ev_1^\ast(E)) + \hwt(F)$ for all $F \in \QC(\mcX)$, and $E \otimes (-)$ preserves $\Cycles(\mcX/\mcB)^{<\infty}$.
\end{cor}
\begin{proof}
	By \Cref{L:highest_weight_condition}, the first claim is equivalent to $E \otimes \QC_\mcS(\mcX)^{<v} \subset \QC_\mcS(\mcX)^{<v+\hwt(E)}$, where $\QC_\mcS(\mcX)^{<v} = \QC_\mcS(\mcX) \cap \QC(\mcX)^{<v}$. By \Cref{L:local_cohomology_formula}, any $F \in \QC_\mcS(\mcX)^{<v}$ has a bounded below convergent filtration whose associated graded pieces are of the form $(\ev_1)_\ast(G)$ for $G \in \QC(\mcS)^{<v}$. Because $\QC_\mcS(\mcX)^{<v+\hwt(E)}$ is closed under filtered colimits, this reduces the first claim to showing that $E \otimes (\ev_1)_\ast(-) \cong (\ev_1)_\ast(\ev_1^\ast(E) \otimes -)$ maps $\QC(\mcS)^{<v}$ to $\QC(\mcX)^{<v+\hwt(\ev_1^\ast(E))}$. This follows from \Cref{L:perf_completion_symmetric_monoidal} and the fact that $(\ev_1)_\ast$ commutes with baric truncation.
	
	Now consider $F \in \Cycles(\mcX/\mcB)^{<\infty}$. By hypotheses, $E \otimes F |_{\mcX \setminus \mcS} \in \Cycles(\mcX \setminus \mcS)$, and we have already shown that $E \otimes F \in \QC(\mcX)^{<\infty}$, so \Cref{L:baric_truncation_preserves_highest_weight_cycles} reduces the claim to showing that $\radj{w}(E \otimes F) \in \Cycles(\mcX/\mcB)^{<\infty}$ for all $w \in \bbZ$. If $\ev_1^\ast(E) \in \QC(\mcS)^{<u}$ for $u>0$, then
	\[
	\radj{w}(E \otimes F) \cong \radj{w}(E \otimes \radj{w-u}(F))
	\]
	Again by \Cref{L:baric_truncation_preserves_highest_weight_cycles}, $\radj{w-u}(F)$ has a finite filtration, the baric filtration, whose associated graded pieces have the form $(\ev_1)_\ast(G)$ for $G \in \Cycles(\mcS/\mcB)^{<\infty}$ by \Cref{L:constant_weight_cycles}. It therefore suffices to show that $\radj{w}(E \otimes (\ev_1)_\ast(-)) \cong (\ev_1)_\ast(\radj{w}(\ev_1^\ast(E) \otimes (-)))$ maps $\Cycles(\mcS/\mcB)^{<\infty}$ to $\Cycles(\mcX/\mcB)^{<\infty}$. This follows from \Cref{L:highest_weight_cycles_module_category}, \Cref{L:baric_truncation_preserves_highest_weight_cycles}, and \Cref{P:functoriality_highest_weight_cycles}.
\end{proof}

\begin{lem}\label{L:highest_weight_cycles_are_cycles}
	$\Cycles(\mcX/\mcB)^{<\infty} \subset \Cycles(\mcX/\mcB)$
\end{lem}
\begin{proof}
    For any $F \in \QC(\mcX)$ and $P \in \Perf(\mcX)$ one has $\ev_1^!(F \otimes P) \cong \ev_1^\ast(P) \otimes \ev_1^!(F)$, so $\Cycles(\mcX/\mcB)^{<\infty}$ is closed under $(-) \otimes P$. To prove the claim it therefore suffices to check that $\rho_\ast(F) \in \DCoh(\mcB)$ for any $F \in \Cycles(\mcX/\mcB)^{<\infty}$. From the exact triangle $R\Gamma_\mcS(F) \to F \to i_\ast(i^\ast(F))$ and the fact that $F|_{\mcX\setminus \mcS} \in \Cycles(\mcX \setminus \mcS)$, this is equivalent to showing $\rho_\ast(R\Gamma_\mcS(F)) \in \DCoh(\mcB)$.

    Because $\rho$ is universally of finite cohomological dimension, $\rho_\ast$ commutes with filtered colimits \cite[Thm.~A.1.5]{MR4560539}.
	Therefore, by \Cref{L:local_cohomology_formula}, $\rho_\ast(R\Gamma_\mcS(F))$ has a convergent filtration whose associated graded pieces are isomorphic to \[\rho_\ast(R\uHom_\mcX((\ev_1)_\ast(\Sym^n(\bbL_{\mcS/\mcX})), F)) \cong (\rho|_\mcS)_\ast( \Sym^n(\bbL_{\mcS/\mcX}^\ast) \otimes \ev_1^!(F)).
    \]
    Because $\bbL_{\mcS/\mcX}^\ast \in \Perf(\mcS)^{<0}$, $\ev_1^!(F) \in \QC(\mcS)^{<v}$ for some $v \in \bbZ$, and $\ev_1^!(F) \in \Cycles(\mcS/\mcB)$, it follows that all of these associated graded pieces lie in $\DCoh(\mcB)$ and that they vanish for $n\geq v$. This implies that $\rho_\ast(R\Gamma_\mcS(F)) \in \DCoh(\mcB)$.
\end{proof}

\subsection{The non-abelian localization theorem}\label{S:main_result}

We continue to work under \Cref{H:lfp}, with $\mcX$ quasi-compact, and we assume that $\mcX$ is equipped with a (necessarily finite) $\Theta$-stratification relative to $\mcB$. It will be convenient to regard the centers $\mcZ = \bigsqcup_i \mcZ_i \subset \Grad(\mcX)$ and the strata $\mcS = \bigsqcup_i \mcS_i \subset \Filt(\mcX)$ as single open substacks. The canonical action of the monoid $\Theta$ on $\mcZ \subset \Grad(\mcX)$ and $\mcS \subset \Filt(\mcX)$ induces a baric structure on both $\QC(\mcZ)$ and $\QC(\mcS)$, as in \Cref{S:baric_stratum}.

We will make frequent use of the canonical morphisms \eqref{E:theta_stratum_arrows}. For the functors $(\ev_1)_\ast$, $\tot_\ast$, $\spfilt_\ast$, and $\spfilt^\ast$, we will use a subscript $i$ to denote the functor for the $i^{\rm{th}}$ stratum, and no subscript to denote the functor from the union of strata. For instance, $(\ev_{1,i})_\ast : \QC(\mcS_i) \to \QC(\mcX)$ and $(\ev_1)_\ast : \QC(\mcS) \to \QC(\mcX)$. The functor $\ev_{1,i}^! : \QC(\mcX) \to \QC(\mcS_i)$ will denote the composition of restriction $\QC(\mcX) \to \QC(\mcX_{\leq i})$ with the $!$-pullback $\QC(\mcX_{\leq i}) \to \QC(\mcS_i)$. We let $\ev_1^! : \QC(\mcX) \to \QC(\mcS) \cong \bigoplus_i \QC(\mcS_i)$ denote the sum of the functors $\ev_{1,i}^!$.

\begin{defn}[Sharp pullback]\label{D:sharp_pullback}
    Let $\bbL^- := \ladj{0}(\tot^\ast(\bbL_{\mcX/\mcR})) \in \QC(\mcZ)$, $\bbL^+ := \radj{1}(\tot^\ast(\bbL_{\mcX/\mcR})) \in \QC(\mcZ)$, and $\bbN_{\mcZ/\mcX} := (\bbL^+ \oplus \bbL^-)^*$. Then we define a functor $\tot^\sharp \cong \bigoplus_i \tot_i^\sharp : \QC(\mcX) \to \QC(\mcZ) \cong \bigoplus_i \QC(\mcZ_i)$ by
\begin{equation}\label{E:def_tot_sharp}
	\tot^{\sharp}(F) := \det(\bbL^+) \otimes \spfilt^\ast(\ev_1^!(F))[\rank \bbL^+].
\end{equation}
In this formula, $\rank(\bbL^+)$ is interpreted as a locally constant function, so the shift might be different on different components of $\mcZ$.
\end{defn}

We will use the standard notation $R\Gamma_{\mcS_i} \mcO_{\mcX}$ and $ \mcO_{\mcX_{\leq i}} \in \QC(\mcX)$ to denote the push forward of these complexes along the open immersion $\mcX_{\leq i} \subset \mcX$, and $R\Gamma_{\mcS_i}(F) \cong F \otimes R\Gamma_{\mcS_i}(\mcO_\mcX)$.

\begin{defn}[Highest weight cycles, multiple strata]\label{D:highest_weight_cycles_multiple}
	For $\mcX$ equipped with a $\Theta$-stratification as above, we define
	\begin{gather*}
        \hwt(F) := \sup \left\{ w \mid \exists d \in \bbZ, i\in \{0,\ldots,N\} \text{ s.t. } \radj{w}(\ev_1^!(\tau_{\leq d}(R\Gamma_{\mcS_i}(F)))) \neq 0 \right\} \in \bbZ \cup \{\pm \infty\}\\
		\QC(\mcX)^{<\infty} := \left\{F \in \QC(\mcX) \mid \hwt(F) < \infty \right\}\\
        \Cycles(\mcX/\mcB)^{<\infty} := 
        \left\{F \in \QC(\mcX)^{<\infty} \mid \ev_1^!(F) \in \Cycles(\mcS/\mcB) \right\}
	\end{gather*}
\end{defn}

As in the previous section, if there is a stratum $\mcS_i$ with a trivial $\Theta$-action, then $\Cycles(\mcS_i)^{<\infty} = \Cycles(\mcS_i)$. Also, if $\mcX$ is eventually coconnective or $F \in \QC(\mcX)_{<\infty}$, then $F \in \QC(\mcX)^{<\infty}$ if and only if $\ev_1^!(F) \in \QC(\mcS)^{<\infty}$ by \Cref{L:highest_weight_condition}.

\begin{lem}\label{L:highest_weight_cycles_are_cycles_multiple_strata}
    $\Cycles(\mcX/\mcB)^{<\infty} \subset \Cycles(\mcX/\mcB)$
\end{lem}

\begin{proof}
    The definitions imply that restriction to the open union of strata $\mcX_{\leq i} \subset \mcX$ preserves both categories $\QC(-)^{<\infty}$ and $\Cycles(-)^{<\infty}$. So the claim follows from recursive application of \Cref{L:highest_weight_cycles_are_cycles}.
\end{proof}

\begin{lem}\label{L:characterize_highest_weight_cycles}
    If $K \in \QC(\mcX)^{<\infty}$, then $K \in \Cycles(\mcX/\mcB)^{<\infty}$ if and only if $\tot^\sharp(K) \in \Cycles(\mcZ/\mcB)$.
\end{lem}
\begin{proof}
    From \Cref{D:highest_weight_cycles_multiple}, this is equivalent to the claim that for $F \in \QC(\mcS)^{<\infty}$, $F \in \Cycles(\mcS/\mcB)$ if and only if $\spfilt^\ast(F) \in \Cycles(\mcZ/\mcB)$. By \Cref{L:constant_weight_category} and the fact that $\spfilt^\ast$ commutes with baric truncation, $\beta^w(F) \cong \gr^\ast(\beta^w(\spfilt^\ast(F)))$ for all $w$, so the equivalence follows from \Cref{P:functoriality_highest_weight_cycles} and \Cref{L:baric_truncation_preserves_highest_weight_cycles} applied to both $\mcZ$ and $\mcS$.
\end{proof}

\subsubsection*{Statement of the theorem}

Recall that by \Cref{L:highest_weight_cycles_module_category}, the ring $K_0(\Perf(\mcZ)^\wedge_\beta)$ acts on $K_0(\Cycles(\mcZ/\mcB))$ via tensor product, and $e(\bbN_{\mcZ/\mcX}) \in K_0(\Perf(\mcZ)^\wedge_\beta)$ is a unit by \Cref{L:euler_class}, where $\bbN_{\mcZ/\mcX}:=\bbL_{\mcZ/\mcX}^\ast[1]$ is the normal complex.

\begin{thm}\label{T:non-abelian_localization}
	Let $\rho : \mcX \to \mcB$ be as in \Cref{H:lfp}, with $\mcX$ quasi-compact and equipped with a $\Theta$-stratification $\mcX = \mcS_0 \sqcup \cdots \sqcup \mcS_N \subset \Filt_\mcR(\mcX)$ relative to $\mcB$ as in \Cref{D:relative_theta_stratification}. Then the functors $\tot_\ast : \QC(\mcZ) \to \QC(\mcX)$ and $\tot^{\sharp} : \QC(\mcX) \to \QC(\mcZ)$ preserve the subcategories of highest weight cycles $\Cycles(-/\mcB)^{<\infty}$, and they induce an isomorphism on $K$-theory
	\[
		K_0(\Cycles(\mcX/\mcB)^{<\infty}) \cong K_0(\Cycles(\mcZ/\mcB)^{<\infty}) \cong \bigoplus_i K_0(\Cycles(\mcZ_i/\mcB))^{<\infty}.
	\]
	For any $E \in \Cycles(\mcX/\mcB)^{<\infty}$ one has an equality in $K_0(\Cycles(\mcX/\mcB)^{<\infty})$
	\begin{equation}\label{E:main_K_theory_statement}
		[E] = \tot_\ast\left(\frac{1}{e(\bbN_{\mcZ/\mcX})} \cdot \tot^{\sharp}([E]) \right)
	\end{equation}
\end{thm}

We will prove this result at the end of the subsection, after explaining it a bit more. The functor $\rho_\ast : \Cycles(\mcX/\mcB)^{<\infty} \to \DCoh(\mcB)$ induces an index homomorphism $\chi(\mcX,-) : K_0(\Cycles(\mcX/\mcB)^{<\infty}) \to K_0(\DCoh(\mcB))$. If $\mcB = \Spec(k)$ for a field, then $K_0(\DCoh(\mcB)) \cong \bbZ$, and this is the usual $K$-theoretic index $E \mapsto \sum_j (-1)^j \dim H_j(R\Gamma(\mcX,E))$. Suppose one has a ``fundamental class'' $\mathbf{O} \in K_0(\Cycles(\mcX/\mcB)^{<\infty})$. Then \Cref{T:non-abelian_localization} combined with the fact that $\tot^\sharp(P \otimes E) \cong \tot^\ast(P) \otimes \tot^\sharp(E)$ for $P \in \Perf(\mcX)$ implies that for any $P \in \Perf(\mcX)$,
\begin{equation}\label{E:perfect_formula}
\chi(\mcX, P \otimes \mathbf{O}) = \sum_i \chi\left(\mcZ_i, \frac{[\tot^\ast(P)]}{e(\bbN_{\mcZ/\mcX})} \cdot \tot^\sharp(\mathbf{O})\right). \end{equation}
This formula is a derived analogue of the $K$-theoretic trace formula. There are two main subtleties in applying the formula at this greater level of generality: identifying a fundamental class $\mathbf{O} \in \Cycles(\mcX)^{<\infty}$ and describing $\tot^\sharp(\mathbf{O})$.

If $\rho : \mcX \to \mcB$ is cohomologically proper, which under our other hypotheses is equivalent to saying that $\rho$ maps $\DCoh(\mcX)$ to $\DCoh(\mcB)$, then $\DCoh(\mcX) \subset \Cycles(\mcX/\mcB)$. If in addition $\mcX$ is eventually coconnective, meaning $H_i(\mcO_\mcX) \cong 0$ for $i \gg 0$, then $\mcO_\mcX \in \DCoh(\mcX) \subset \Cycles(\mcX/\mcB)$. We can therefore regard $\mcO_\mcX$ as a fundamental class in $K_0(\Cycles(\mcX/\mcB))$.

\begin{ex}
	If $\mcX$ is an eventually coconnective algebraic derived stack of finite presentation that admits a proper relative good moduli space $X$ over $\mcB$, then it is cohomologically proper, and hence $\mcO_\mcX \in \Cycles(\mcX/\mcB)$. This is because the pushforward $q_\ast : \QC(\mcX) \to \QC(X)$ preserves $\DCoh(-)$ and $\DCoh(X) \subset \Cycles(X/\mcB)$. 
\end{ex}

We have a satisfactory calculation of $\tot^\sharp(\mathbf{O})$ when $\pi_\mcX$ is quasi-smooth, meaning $\bbL_{\mcX/\mcR}$ has Tor amplitude in $[-1,1]$. We note that if $\mcX \to \mcR$ is quasi-smooth and $\mcR$ is eventually coconnective, then so is $\mcX$.
\begin{cor}[Quasi-smooth localization formula]\label{T:quasi-smooth_localization}
In the context of \Cref{T:non-abelian_localization}, if $\pi_\mcX : \mcX \to \mcR$ is quasi-smooth and $\mcR$ is eventually coconnective and almost finitely presented over a field $k$ of characteristic $0$, then $\mcO_{\mcX} \in \QC(\mcX)^{<\infty}$ and $\tot^\sharp(\mcO_\mcX) \cong \mcO_\mcZ$. Suppose in addition that $\mcZ_i$ is cohomologically proper over $\mcB$ for all $i$, which is automatically true if $\mcX$ admits a proper relative good moduli space over $\mcB$. Then $\Perf(\mcX) \subset \Cycles(\mcX/\mcB)^{<\infty}$, and for any $P \in \Perf(\mcX)$,
	\begin{equation}\label{E:quasi-smooth_localization_formula}
		\chi(\mcX, P) = \sum_i \chi\left(\mcZ_i, e(\bbN_{\mcZ_i/\mcX})^{-1} \cdot [\tot_i^\ast(P)]\right),
	\end{equation}
	where $e(N_{\mcZ_i/\mcX})^{-1}$ is the inverse of the Euler class in $K_0(\Perf(\mcZ_i)^\wedge_\beta)$. (See \Cref{L:highest_weight_cycles_module_category}.)
\end{cor}
\begin{proof}
    We will make use of the formalism of Ind-coherent sheaves on locally almost finitely presented $k$-stacks developed in \cite{MR3701352, Arinkin2015SingularConjecture,MR3136100}. Because $\mcR$ is eventually coconnective and $\pi_\mcX$ is quasi-smooth, $\mcX$ is eventually coconnective as well. Because $\ev_1$ is a closed immersion, the shriek-pullback functor on $\ICoh$, $\ev_1^{IC,!} : \ICoh(\mcX) \to \ICoh(\mcS)$ agrees with $\ev_1^! : \QC(\mcX) \to \QC(\mcS)$ on the subcategory $\QC(\mcX)_{<\infty} \cong \ICoh(\mcX)_{<\infty}$ of homologically bounded above complexes. We will therefore compute $\ev_1^{IC,!}(\mcO_\mcX)$.

    We regard $\mcO_\mcR \in \DCoh(\mcR) \subset \ICoh(\mcR)$, and we define the relative canonical complex $\omega_{\mcX/\mcR} := \pi_\mcX^!(\mcO_\mcR)$. A quasi-smooth morphism is Gorenstein, meaning $\omega_{\mcX/\mcR}$ is a shift of an invertible sheaf and thus lies in the subcategory $\Perf(\mcX) \subset \DCoh(\mcX) \subset \ICoh(\mcX)$. We have
	\[\ev_1^{IC,!}(\mcO_\mcX) \cong \ev_1^{IC,!}(\omega_{\mcX/\mcR} \otimes \omega_{\mcX/\mcR}^{-1}) \cong \omega_{\mcS/\mcR} \otimes \ev_1^\ast(\omega_{\mcX/\mcR}^{-1}),\]
    where we have identified $\ev_1^{IC,!}(\omega_{\mcS/\mcR}) \cong (\pi_{\mcX}|_{\mcS})^{IC,!}(\mcO_\mcR) \cong \omega_{\mcS/\mcR}$ and the fact that shriek pullback in $\ICoh$ is compatible with the action of $\QC$ on $\ICoh$ by tensor product. $\mcS \to \mcR$ is also quasi-smooth by \Cref{L:stratum_cotangent}. Forthcoming work of Adeel Kahn will show that by applying derived deformation to the normal bundle \cite{normal_cone} to the quasi-smooth morphisms $\mcX \to \mcR$ and $\mcS \to \mcR$ one can construct an isomorphism (See also \cite[App.~B]{Halpern-Leistner2021DerivedConjecture} for a proof for quotient stacks and a slightly weaker statement in general.)
	\[
	\omega_{\mcX/\mcR} \cong \det(\bbL_{\mcX/\mcR})[\rank \bbL_{\mcX/\mcR}] \text{ and } \omega_{\mcS/\mcR} \cong \det(\bbL_{\mcS/\mcR})[\rank \bbL_{\mcS/\mcR}].
	\]
    Applying $\det$ to the exact triangle $\ev_1^\ast(\bbL_{\mcX/\mcR}) \to \bbL_{\mcS/\mcR} \to \bbL_{\mcS/\mcX}$, we obtain an isomorphism
	\begin{align*}
	\ev_1^!(\mcO_\mcX) &\cong \det(\bbL_{\mcS/\mcR}) \otimes \ev_1^\ast(\det(\bbL_{\mcX/\mcR})^{-1}) [\rank \bbL_{\mcS/\mcR} - \rank \bbL_{\mcX/\mcR}] \\
	&\cong \det(\bbL_{\mcS/\mcX})[\rank \bbL_{\mcS/\mcX}].
	\end{align*}
	It follows from \Cref{D:sharp_pullback} that
	\[
	\tot^\sharp(\mcO_\mcX) \cong \det(\bbL^+) \otimes \spfilt^\ast(\det(\bbL_{\mcS/\mcX}))[\rank \bbL_{\mcS/\mcX}+\rank \bbL^+].
	\]
	To complete the proof that $\tot^\sharp(\mcO_\mcX) \cong \mcO_\mcZ$, we use that $\spfilt^\ast(\bbL_{\mcS/\mcX}) \cong \bbL^+[1]$ by \Cref{L:stratum_cotangent}, which implies that $\spfilt^\ast(\det(\bbL_{\mcS/\mcX})) \cong \det(\bbL^+)^{-1}$ and $\rank \bbL^+ = - \rank \bbL_{\mcS/\mcX}$.
    
    To show that $\Perf(\mcX) \subset \Cycles(\mcX/\mcB)^{<\infty}$ if each center $\mcZ_i$ is cohomologically proper, it suffices to show that $\mcO_\mcX \in \Cycles(\mcX/\mcB)^{<\infty}$ in this case. Our computation of $\ev_1^!(\mcO_\mcX)$ shows that $\mcO_{\mcX} \in \QC(\mcX)^{<\infty}$, so the claim follows from \Cref{L:characterize_highest_weight_cycles}. The formula in the statement is just \Cref{E:perfect_formula} combined with the fact that $\tot^\sharp(\mcO_\mcX) \cong \mcO_\mcZ$.
\end{proof}

\subsubsection*{Completing the proof of \Cref{T:non-abelian_localization}}

The diagram $\mcO_\mcX = \mcO_{\mcX_{\leq N}} \to \mcO_{\mcX_{\leq n-1}} \to \cdots \to \mcO_{\mcX_{\leq 0}}$ constitutes a finite filtration of $\mcO_\mcX$ whose associated graded complexes are $R\Gamma_{\mcS_i} \mcO_\mcX$ for $i=0,\ldots,N$. (Recall that $\mcS_0 = \mcX_{\leq 0}$ by convention, so $R\Gamma_{\mcS_0} \mcO_{\mcX} =\mcO_{\mcX_{\leq 0}}$.) We call this the local cohomology filtration of $\mcO_\mcX$.

\begin{lem}\label{L:local_cohomology_sod}
	The filtration of the identity functor $F \mapsto \mcO_\mcX \otimes F$ coming from the local cohomology filtration of $\mcO_\mcX$ is the filtration associated to a semiorthogonal decomposition
	\[
	\Cycles(\mcX/\mcB)^{<\infty} = \langle \mcA_0,\ldots,\mcA_N \rangle,
	\]
	where $\mcA_i = R\Gamma_{\mcS_i}(\mcO_\mcX) \otimes \Cycles(\mcX/\mcB)^{<\infty}$ is the image under the pushforward $\QC(\mcX_{\leq i}) \to \QC(\mcX)$ of the full subcategory of complexes in $\Cycles(\mcX_{\leq i})^{<\infty}$ that are set theoretically supported on $\mcS_{i}$.
\end{lem}
\begin{proof}
For any $E,F \in \QC(\mcX)$ and $i<j$, one has $R\Hom_\mcX(R\Gamma_{\mcS_j}(\mcO_\mcX) \otimes E, R\Gamma_{\mcS_i}(\mcO_\mcX) \otimes F) \cong 0$ because the first is set-theoretically supported on $\mcS_j$ and the second is the pushforward of a complex along the open immersion $\mcX_{\leq j-1} = \mcX_{\leq j} \setminus \mcS_j \subset \mcX_{\leq j}$. The facts that $R\Gamma_{\mcS_i}(\mcO_\mcX) \otimes R\Gamma_{\mcS_j}(\mcO_\mcX) \cong \mcO_\mcX^{\delta_{i,j}} \otimes R\Gamma_{\mcS_i}(\mcO_\mcX)$ and $\ev_{1,i}^!(R\Gamma_{\mcS_j}(\mcO_\mcX) \otimes E) \cong \mcO_\mcS^{\delta_{i,j}} \ev_{1,i}^!(E)$ imply that $R\Gamma_{\mcS_i}(\mcO_\mcX) \otimes (-)$ preserves $\Cycles(\mcX/\mcB)^{<\infty}$. The claim follows.
\end{proof}

\begin{proof}[Proof of \Cref{T:non-abelian_localization}]
	\Cref{P:functoriality_highest_weight_cycles} implies that $(\ev_{1,i})_\ast \circ (\spfilt_{i})_\ast : \QC(\mcZ_i) \to \QC(\mcX)$ maps $\Cycles(\mcZ_i)^{<\infty}$ to $\mcA_i$. In the other direction, $\spfilt_i^\ast \circ \ev_{1,i}^! : \Cycles(\mcX/\mcB)^{<\infty} \to \QC(\mcZ_i)$ evidently factors through the projection $R\Gamma_{\mcS_i} \mcO_\mcX \otimes (-) : \Cycles(\mcX/\mcB)^{<\infty} \to \mcA_i$ and has image in $\Cycles(\mcZ_i)^{<\infty}$ by \Cref{P:functoriality_highest_weight_cycles} applied to the stratum $\mcS_i \hookrightarrow \mcX_{\leq i}$. This implies the same for $\tot^\sharp$. \Cref{L:local_cohomology_sod} implies that $\bigoplus_i \mcA_i \to \Cycles(\mcX/\mcB)^{<\infty}$ induces an isomorphism on $K$-theory, so it suffices to show that $\Cycles(\mcZ_i)^{\infty} \to \mcA_i$ is an isomorphism on $K_0$ for each $i$ individually. We will therefore drop the subscripts and just consider the case where $\mcS \hookrightarrow \mcX$ is a closed $\Theta$-stratum.

	We first claim that $(\ev_1)_\ast : \Cycles(\mcS/\mcB)^{<\infty} \to \mcA$, $\gr^\ast : \Cycles(\mcZ/\mcB)^{<\infty} \to \Cycles(\mcS/\mcB)^{<\infty}$, and $\spfilt^\ast : \Cycles(\mcS/\mcB)^{<\infty} \to \Cycles(\mcZ/\mcB)^{<\infty}$ induce isomorphisms on $K_0$, with the last two being inverse to each other. \Cref{L:baric_truncation_preserves_highest_weight_cycles} and \Cref{L:decompositions_of_baric_completion} imply that the baric structure on $\QC(\mcX)$ restricts to a right-complete baric structure on $\mcA = (\Cycles(\mcX/\mcB)^{<\infty})^{\rm{nil}} \subset \Cycles(\mcX/\mcB)^{\infty}$, which is the full subcategory of objects supported set theoretically on $\mcS$ by \Cref{L:first_local_cohomology_formula}. Regarding $\mcZ$ and $\mcS$ tautologically as $\Theta$-strata in themselves, the same argument shows that the baric structures on $\Cycles(\mcZ/\mcB)^{<\infty}$ and $\Cycles(\mcS/\mcB)^{<\infty}$ are right-complete. The functors $\gr^\ast$ and $(\ev_1)_\ast$ commute with baric truncation, so \Cref{L:baric_splitting_K_theory} implies that the claim for $(\ev_1)_\ast$ and $\gr^\ast$ is equivalent to showing the same for $\Cycles(\mcS/\mcB)^w \to \mcA^w \cong \QC(\mcX)^w$ and $\Cycles(\mcZ/\mcB)^w \to \Cycles(\mcS/\mcB)^w$. These last two functors are isomorphisms by \Cref{L:constant_weight_cycles}.

    To show that $\spfilt_\ast : \Cycles(\mcZ/\mcB)^{<\infty} \to \Cycles(\mcS/\mcB)^{<\infty}$ induces an isomorphism on $K_0$, it now suffices to show that the composition $\spfilt^\ast \circ \spfilt_\ast$ induces an isomorphism on $K_0$. We claim that for any $F \in \Cycles(\mcS/\mcB)^{<\infty}$, $[\spfilt^\ast \circ \spfilt_\ast(F)] = [\Sym(\bbL^-[1]) \otimes F]$. By \Cref{L:baric_splitting_K_theory}, this is equivalent to showing that $[\radj{v}(\spfilt^\ast \circ \spfilt_\ast(F))] = [\radj{v}(\Sym(\bbL^-[1]) \otimes F)]$ for all $v \in \bbZ$, which is an immediate consequence of \Cref{T:push_pull_filtration_stratum} and the fact that $\bbL_{\mcZ/\mcS} \cong \bbL^-[1] \in \Perf(\mcZ)^{<0}$ by \Cref{L:stratum_cotangent}. Now $e((\bbL^-)^\ast) = [\Sym(\bbL^{-}[1])] \in K_0(\Perf(\mcS)^\wedge_\beta)$ is a unit by \Cref{L:euler_class}, so $\spfilt_\ast \circ \spfilt^\ast = e((\bbL^-)^\ast) \cdot (-)$ induces an isomorphism on $K_0$, and therefore so does $\spfilt_\ast$.

	Finally, we describe the isomorphism $\tot_\ast = (\ev_1)_\ast \circ \spfilt_\ast$ on $K_0$ more explicitly. The previous paragraph shows that
	\[[E] = \left[\spfilt_\ast\left(\frac{1}{e((\bbL^-)^\ast)} \cdot \spfilt^\ast(E)\right)\right]\]
	in $K_0(\Cycles(\mcS/\mcB)^{<\infty})$. On the other hand, by \Cref{L:local_cohomology_formula}, any $F \in \mcA \subset \QC(\mcX)^{<\infty}$ has a bounded below convergent filtration $F \cong \colim_n R\uHom(\mcO_{\mcS^{(n)}}, F)$ with graded pieces $G_n \cong (\ev_1)_\ast(\Sym^n(\bbL_{\mcS/\mcX}^\ast [1]) \otimes \ev_1^!(F))$. One has $\bbL_{\mcS/\mcX} \cong \radj{1}(\ev_1^\ast(\bbL_{\mcX/\mcR})[1]) \in \Perf(\mcS)^{\geq 1}$ by \Cref{L:stratum_cotangent}. Using \Cref{L:K_theory_graded} as in the last paragraph, one can then show that
	\[[F] = \left[(\ev_1)_\ast\left(\Sym\left(\left(\radj{1}(\ev_1^\ast(\bbL_{\mcX/\mcR}))\right)^\ast\right) \otimes \ev_1^!(F)\right)\right]\]
	in $K_0(\mcA)$. Observing that $\spfilt^\ast(\radj{1}(\ev_1^\ast \bbL_{\mcX/\mcR})) \cong \bbL^+$, we can combine this with the previous formula to obtain
	\[
		[F] = \left[\tot_\ast\left(\frac{\Sym((\bbL^+)^\ast)}{e((\bbL^-)^\ast)} \cdot \spfilt^\ast(\ev_1^!(F))\right)\right].
	\]
	The final formula in the statement of \Cref{T:non-abelian_localization} follows from the observation that $\frac{1}{e((\bbL^+)^\ast)} = e((\bbL^+)^\ast[-1]) = \Sym((\bbL^+)^\ast) \otimes \det(\bbL^+)^\ast[-\rank(\bbL^+)]$. In order to arrive at a more symmetric formula, we cancel the last two terms by adding a factor of $\det(\bbL^+)[\rank(\bbL^+)]$ to $\tot^\sharp$.
\end{proof}

\subsection{Completion of \texorpdfstring{$\Perf(\mcX)$}{Perf(X)}}

Suppose $\rho : \mcX \to \mcB$ is as in \Cref{H:lfp}, with $\mcX$ quasi-compact. Suppose $\mcX$ is equipped with a $\Theta$-stratification relative to $\mcB$ as in \Cref{D:relative_theta_stratification}. We will consider the following category, which behaves like a categorification of a completion of cohomology, rather than homology:

\begin{defn}
	The \emph{completion} of $\Perf(\mcX)$ is the full subcategory
	\[
	\Perf(\mcX)^\wedge_\beta := \left\{F \in \QC(\mcX) \left| \ev_1^\ast(F) \in \Perf(\mcS)_\beta^\wedge \right. \right\}.
	\]
\end{defn}

In the degenerate case where $\mcS \cong \mcX$, i.e., $\mcX$ is a single $\Theta$-stratum, this definition is equivalent to the definition of $\Perf(\mcS)^\wedge_\beta$ studied in \Cref{S:baric_completion_perf}.

\begin{ex}
Unlike the category of cycles, the category $\Perf(\mcX)^\wedge_\beta$ depends on the choice of $\Theta$-stratification. For instance consider $\mcX = \bbA^1 / \bbG_m$, where the action arises by equipping the polynomial ring $k[x]$ with a grading in which $x$ has degree $-1$, where $k$ is a field. If we regard $\mcX$ as a single $\Theta$-stratum, then we have seen in \Cref{EX:perf_completion_Z} (using \Cref{L:perf_completion_functoriality_2}) that $\Perf(\mcX)^\wedge_\beta$ is equivalent to the category of complexes of graded $k[x]$-modules $F^\bullet$ such that the weights of $H^\ast(F^\bullet)$ are $<v$ for some $v < \infty$, and such that $H^\ast(F^\bullet \otimes_{k[x]} k)$ is finite dimensional in each weight. On the other hand, if we consider the $\Theta$-stratification $\mcX \cong \{0\}/\bbG_m \sqcup (\bbA^1 \setminus 0) / \bbG_m$, then $\Perf(\mcX)^\wedge$ is equivalent to the category of graded $k[x]$-modules $F^\bullet$ such that $H^\ast(F \otimes_{k[x]} k)$ is finite dimensional in each weight and only nonzero in weights $\geq w$ for some $w \in \bbZ$.
\end{ex}

\begin{prop}\label{P:perf_completion}
The category $\Perf(\mcX)^\wedge_\beta$ has the following properties:
\begin{enumerate}
	\item If $\mcU \subset \mcX$ is an open union of $\Theta$-strata, then the functors $j_\ast : \QC(\mcU) \to \QC(\mcX)$, $j^\ast : \QC(\mcX) \to \QC(\mcU)$, and $R\underline{\Gamma}_{\mcX \setminus \mcU}(-)$ preserve $\Perf(-)^\wedge_\beta$.
	\item For $E \in \Perf(\mcX)^\wedge_\beta$, $E \otimes (-)$ preserves $\Perf(\mcX)^\wedge_\beta$, $\QC(\mcX)^{<\infty}$, and $\Cycles(\mcX/\mcB)^{<\infty}$. In particular, $\Perf(\mcX)^\wedge_\beta$ is symmetric monoidal, and $\rho_\ast$ defines a pairing $\Perf(\mcX)^\wedge_\beta \otimes \Cycles(\mcX/\mcB)^{<\infty} \to \DCoh(\mcB)$.
\end{enumerate}
\end{prop}
\begin{proof}
	The first condition is automatic from the definition, because for the center of each stratum, i.e., each component $\mcS_i \subset \mcS$, the functors $j_\ast$, $j^\ast$, and $R\underline{\Gamma}$ either preserve $(\ev_{1,i})^\ast(F)$ or set it to $0$. The fact that for $E \in \Perf(\mcX)^\wedge_\beta$, $E \otimes (-)$ preserves $\Perf(\mcX)^\wedge_\beta$ follows from the definition and \Cref{L:perf_completion_symmetric_monoidal}.

	To show that $E \otimes (-)$ preserves $\Cycles(\mcX/\mcB)^{<\infty}$ and $\QC(\mcX)^{<\infty}$, \Cref{L:local_cohomology_sod} allows one to reduce the claim to the situation where $\mcS$ is a single closed $\Theta$-stratum, and for objects in $\Cycles(\mcX/\mcB)^{<\infty}$ or $\QC(\mcX)^{<\infty}$ that are set theoretically supported on $\mcS$. This case was handled in \Cref{C:completion_tensor_property}.
\end{proof}

As a result of \Cref{P:perf_completion}, $K^0(\Perf(\mcX)^\wedge_\beta)$ is a ring, and $1 = [\mcO_\mcX] = \sum_\alpha [R\Gamma_{\mcS_\alpha} \mcO_\mcX]$ is a decomposition of the identity into a sum of mutually orthogonal idempotent classes, which is a version of formulation \ref{I:localization_cohomology} of the localization theorem in the introduction. To relate this to the centers of the strata, we will need a stronger hypothesis on the $\Theta$-stratification.

We will say the $\Theta$-stratification is \emph{regular} if the inclusion $\ev_1 : \mcS_{\alpha} \hookrightarrow \mcX \setminus \bigcup_{\beta >\alpha} \mcS_{\beta}$ is a regular embedding for all $\alpha$. This is automatic when $\mcX$ is smooth. The regular embedding condition implies that $\ev_1^!(-) \cong \ev_1^\ast(-) \otimes \det(\bbL_{\mcS/\mcX})[\rank(\bbL_{\mcS/\mcX})]$, which simplifies the description of many of the categories considered above.

\begin{prop}\label{P:perf_completion_regular}
If the $\Theta$-stratification of $\mcX$ is regular, then the functors $\tot^\ast$, $\tot_\ast$, $\spfilt^\ast$, $\spfilt_\ast$, $\ev_1^\ast$ and $(\ev_1)_\ast$ preserve the categories $\Perf(-)^\wedge_\beta$ and induce isomorphisms
	\[K_0(\Perf(\mcX)^\wedge_\beta) \overset{\ev_1^\ast}{\cong} K_0(\Perf(\mcS)^\wedge_\beta) \overset{\spfilt^\ast}{\cong} K_0(\Perf(\mcZ)^\wedge_\beta).
	\]
	For any $[E] \in K_0(\Perf(\mcX)^\wedge_\beta)$, $[E] = \tot_\ast(e(\bbN_{\mcZ/\mcX})^{-1} \cdot \tot^\ast([E])).$
\end{prop}

\begin{proof}
The functoriality for $\spfilt^\ast$ and $\spfilt_\ast$ is established in \Cref{L:perf_completion_functoriality} and \Cref{L:perf_completion_functoriality_2}, and  $\spfilt^\ast$ induces an isomorphism $K_0(\Perf(\mcS)^\wedge_\beta) \cong K_0(\Perf(\mcZ)^\wedge_\beta)$ by \Cref{L:K_theory_perf_completion}.
	
	The functor $\ev_1^\ast$ preserves $\Perf(-)^\wedge_\beta$ by definition. To show that $(\ev_1)_\ast$ preserves $\Perf(-)^\wedge_\beta$, we show that $\ev_1^\ast((\ev_1)_\ast(-))$ preserves $\Perf(\mcS)^\wedge_\beta$. This follows from the fact that, because $\ev_1$ is a regular closed immersion, $\ev_1^\ast((\ev_1)_\ast(E))$ admits a finite filtration with $\gr(\ev_1^\ast((\ev_1)_\ast(E))) \cong E \otimes \Sym(\bbL_{\mcS/\mcX})$, where $\Sym(\bbL_{\mcS/\mcX}) \in \Perf(\mcS)$ again because $\mcS \hookrightarrow \mcX$ is regular.
	
	The same argument as in the proof of \Cref{T:non-abelian_localization}, using the complete baric structures and \Cref{L:baric_splitting_K_theory}, shows that $(\ev_1)_\ast$ induces an isomorphism $K_0(\Perf(\mcS)^\wedge_\beta) \cong K_0(\Perf_\mcS(\mcX)^\wedge_\beta)$, where $\Perf_\mcS(\mcX)^\wedge_\beta := \Perf(\mcX)^\wedge_\beta \cap \QC_\mcS(\mcX)$. On the level of $K$-theory, we have already seen that $\ev_1^\ast((\ev_1)_\ast([E])) = [E] \otimes e(\bbL_{\mcS/\mcX}^\ast[1])$. $e(F)$ is a unit in $K_0(\Perf(\mcS)^\wedge_\beta)$ when it is defined, so $\ev_1^\ast((\ev_1)_\ast(-))$ is an isomorphism on $K_0$, which implies the same for $\ev_1^\ast$. Likewise, the same argument as in the proof of \Cref{T:non-abelian_localization} implies that $\spfilt^\ast(\spfilt_\ast([E])) \cong [E] \cdot e(\bbL_{\mcZ/\mcS}^*[1])$. Combining these calculations shows that $\tot^\ast(\tot_\ast([E])) = [E] \cdot e(\bbN_{\mcZ/\mcX})$, which establishes the final claim of the proposition.
\end{proof}

\begin{ex}
    When $\mcO_\mcX \in \Cycles(\mcX/\mcB)$, such as in the context of \Cref{T:quasi-smooth_localization}, \Cref{P:perf_completion}(2) implies that $\Perf(\mcX)^\wedge_\beta \subseteq \Cycles(\mcX/\mcB)$. In particular, there is an ``Euler characteristic'' homomorphism $\chi : K_0(\Perf(\mcX)^\wedge_\beta) \to K_0(\DCoh(\mcB))$ that extends the Euler characteristic on $K_0(\Perf(\mcX))$. This shows that even though $\Perf(\mcX)^\wedge_\beta$ contains certain infinite sums, the $K$-theory of this category is non-zero.
\end{ex}

\subsection{Generalization to infinite \texorpdfstring{$\Theta$}{Theta}-stratifications} \label{S:infinite_stratifications}

In applications, one sometimes wishes to apply non-abelian localization to a stack $\mcX$ that is locally of finite presentation but not quasi-compact, such as the stack of Higgs bundles on a smooth projective curve \cite{halpernleistner2016equivariantverlindeformulamoduli}, or the stack of sheaves on a surface (see \Cref{S:1D_sheaves} below). One can of course imitate \Cref{D:highest_weight_cycles_multiple}, but \Cref{L:highest_weight_cycles_are_cycles_multiple_strata} would fail if there are infinitely many strata.

The solution is to introduce subcategories of admissible objects, following \cite{halpernleistner2016equivariantverlindeformulamoduli, Teleman2009TheCurve}. We bake these admissibility conditions into our definitions in the context of infinitely many strata.

\begin{defn}\label{D:perf_completion}
Suppose $\rho : \mcX \to \mcB$ is a morphism as in \Cref{H:lfp}, and suppose $\mcX$ has a (potentially infinite) $\Theta$-stratification $\mcX = \bigcup_{\alpha \in I} \mcS_\alpha$ relative to $\mcB$ such that for any $\alpha \in I$, the open substack $\mcX_{\leq \alpha} := \bigcup_{\gamma \leq \alpha} \mcS_\gamma \subset \mcX$ is quasi-compact. Then we define
\[
\Cycles(\mcX/\mcB)^{<\infty} := \left\{F \in \QC(\mcX) \left| \begin{array}{l} F|_{\mcX_{\leq \alpha}} \in \Cycles(\mcX_{\leq \alpha}/\mcB)^{<\infty} \text{ for all }\alpha, \text{ and}\\ \hwt(R
\Gamma_{\mcS_\alpha}(F)) \geq 0 \text{ for only finitely many }\alpha \end{array} \right.\right\} 
\]
\[
\Perf(\mcX)^\wedge_\beta := \left\{F \in \QC(\mcX) \left| \begin{array}{l} \ev_{1,\alpha}^*(F) \in \Perf(\mcS_\alpha)^\wedge_\beta \text{ for all }\alpha, \text{ and}\\ \hwt(\ev_{1,\alpha}^*(F)) > 0 \text{ for only finitely many }\alpha \end{array} \right.\right\} 
\]
\end{defn}

We note that the quasi-compactness of all $\mcX_{\leq \alpha}$ is equivalent to saying that the center of every stratum is quasi-compact and $\{\gamma \in I| \gamma \leq \alpha\}$ is finite for all $\alpha$.

\begin{lem}\label{lem:perf_completion_infinite_strata}
$\Perf(\mcX)^\wedge_\beta \subset \QC(\mcX)$ is a symmetric monoidal subcategory, it preserves $\Cycles(\mcX/\mcB)^{<\infty}$ under $\otimes$, and $\rho_\ast : \QC(\mcX) \to \QC(\mcB)$ maps $\Cycles(\mcX/\mcB)^{<\infty}$ to $\DCoh(\mcB)$. In particular $\Cycles(\mcX/\mcB)^{<\infty} \subset \Cycles(\mcX/\mcB)$.
\end{lem}
\begin{proof}
The only part of the proposition that does not follow from \Cref{T:non-abelian_localization} and the definitions is the last claim, that $\rho_\ast$ maps $\Cycles(\mcX/\mcB)^{<\infty}$ to $\DCoh(\mcB)$. If $\mcS \hookrightarrow \mcX$ is a closed $\Theta$-stratum and $F \in \QC(\mcX)^{<0}$, then applying \cite[Prop.~2.1.4]{Halpern-Leistner2021DerivedConjecture} locally over $\mcB$ shows that the restriction homomorphism $\rho_\ast(F) \to (\rho|_{\mcX \setminus \mcS})_\ast(F|_{\mcX \setminus \mcS})$ is an isomorphism. It follows that the limit
\[
\rho_\ast(F) \cong \lim_{\alpha} (\rho|_{\mcX_{\leq \alpha}})_\ast(F|_{\mcX_{\leq \alpha}}) \cong (\rho|_{\mcX_{\leq \alpha_0}})_\ast(F|_{\mcX_{\leq \alpha_0}})
\]
for some sufficiently large $\alpha_0$. This lies in $\DCoh(\mcB)$ by \Cref{L:highest_weight_cycles_are_cycles_multiple_strata}, applied to the finite $\Theta$-stratification of $\mcX_{\leq \alpha_0}$.
\end{proof}

\begin{lem}\label{L:perfect_highest_weight_cycles_quasi_smooth}
    In the context of \Cref{D:perf_completion}, suppose $\mcR$ is eventually coconnective and almost of finite presentation over a field $k$ of characteristic $0$, and that $\pi_\mcX : \mcX \to \mcR$ is quasi-smooth. If the centers of the $\Theta$-stratification are all cohomologically proper over $\mcB$, then $F \in \Perf(\mcX)$ lies in $\Cycles(\mcX/\mcB)^{<\infty}$ if and only if $\hwt(\ev_{1,\alpha}^\ast(F)) < \eta_\alpha$ for all but finitely many $\alpha$, where \begin{equation}\label{E:define_eta}
    \eta_\alpha := \rm{wt}(\det(\radj{1}(\bbL_{\mcX/\mcR}|_{\mcZ_\alpha}))).
    \end{equation}
\end{lem}
\begin{proof}
    As in the proof of \Cref{T:quasi-smooth_localization}, $\ev_{1,\alpha}^!(\mcO_\mcX) \cong \omega_{\mcS_\alpha/\mcR} \otimes \ev_{1,\alpha}^\ast(\omega_{\mcX/\mcR}^{-1}) \cong \det(\bbL_{\mcS_\alpha / \mcX})[\rank(\bbL_{\mcS_\alpha/\mcX})]$ is a graded invertible sheaf concentrated in weight $-\eta_\alpha$ for each $\alpha$. The claim then follows from: 1) the formula $\ev_{1,\alpha}^!(F) \cong \ev_{1,\alpha}^\ast(F) \otimes \ev_{1,\alpha}^!(\mcO_{\mcX})$, which does not hold in general but holds for $F \in \Perf(\mcX)$ because $F$ is dualizable; and 2) the hypothesis that the center of each stratum is cohomologically proper over $\mcB$, which implies that $\Perf(\mcS_\alpha) \subset \Cycles(\mcS_\alpha/\mcB)$ for all $\alpha$.
\end{proof}

\begin{ex}\label{E:fundamental_class}
	In the context of \Cref{L:perfect_highest_weight_cycles_quasi_smooth}, if $\pi_\mcX$ is smooth, then $\radj{1}(\bbL_{\mcX/\mcR}|_{\mcZ})$ is a locally free sheaf, so $\eta_\alpha\geq 0$ for all $\alpha$ and positive for any non-open stratum. It follows that $\mathbf{O} := \mcO_{\mcX} \in \Cycles(\mcX/\mcB)^{<\infty}$.
\end{ex}
	
When $\mcX$ is not smooth, $\eta_\alpha$ can be positive or negative, so it is not automatic that $\mcO_\mcX \in \Cycles(\mcX/\mcB)^{<\infty}$. In some examples, such as the stack of Higgs bundles \cite{halpernleistner2016equivariantverlindeformulamoduli} or \Cref{S:1D_sheaves}, there is a line bundle $\mcL$ on $\mcX$ such that $\rm{wt}(\mcL|_{\mcZ_\alpha})$ grows negative quickly enough with $\alpha$ that $\rm{wt}(\mcL_\alpha) < \eta_\alpha$ for all but finitely many $\alpha$. In this case $\mathbf{O} := \mcL$ plays the role of a fundamental class.

Finally, when $\rho$ is not quasi-compact, several fundamental properties of the pushforward $\rho_\ast : \QC(\mcX) \to \QC(\mcB)$ fail, but they do hold for categories of highest weight complexes. We omit proofs for the following three lemmas, because their proof uses the same idea as the proof of \Cref{lem:perf_completion_infinite_strata}: the vanishing theorem of \cite[Prop.~2.1.4]{Halpern-Leistner2021DerivedConjecture} reduces the claim to the analogous claim for $\mcX_{\leq \alpha}$ for some $\alpha$, in which case the morphism is quasi-compact and the claim is known.

\begin{lem}
    Let $\{E_i\}_{i \in D}$ be a filtered diagram of objects $E_i \in \QC(\mcX)$ such that $\exists i_0 \in D$ and $\alpha_0 \in I$ such that $\hwt(\ev_{1,\alpha}^!(E_i)) < 0$ for all $\alpha > \alpha_0$ and all $i$ admitting a map $i_0 \to i$. Then the canonical homomorphism is an isomorphism
    \[
    \colim_{i \in D} \rho_\ast(E_i) \to \rho_\ast(\colim_{i \in D} E_i).
    \]
\end{lem}

\begin{lem}
    Suppose $E \in \QC(\mcX)$ is such that $\hwt(\ev_{1,\alpha}^!(E))<0$ for all but finitely many $\alpha$ and $F \in \QC(\mcB)$, then the canonical homomorphism is an isomorphism
    \[
    \rho_\ast(E) \otimes F \to \rho_\ast(\rho^\ast(F) \otimes E).
    \]
\end{lem}
\begin{lem}
    Suppose $E \in \QC(\mcX)$ is such that $\hwt(\ev_{1,\alpha}^!(E))<0$ for all but finitely many $\alpha$, and consider a cartesian square of algebraic derived stacks, with $\mcB'$ noetherian,
    \[
    \xymatrix{\mcX' \ar[r]^{\phi'} \ar[d]^{\rho'} & \mcX \ar[d]^\rho \\ \mcB' \ar[r]^{\phi} & \mcB}.
    \]
    Then the canonical homomorphism is an isomorphism $\phi^\ast(\rho_\ast(E)) \to (\rho')_\ast((\phi')^\ast(E))$.
\end{lem}
\subsubsection*{A variant on admissibility}

In practice, it is useful to have a notion of admissibility \emph{relative} to a choice of fundamental class:

\begin{defn}
	In the context of \Cref{D:perf_completion}, for a fixed class $\mathbf{O} \in \Cycles(\mcX/\mcB)^{<\infty}$, we define the $\mathbf{O}$-\emph{admissible perfect complexes} to be the full subcategory $\Perf(\mcX)^{\mathbf{O}\rm{-adm}} \subset \Perf(\mcX)$ of perfect complexes $E$ such that $E^{\otimes m} \otimes \mathbf{O} \in \Cycles(\mcX/\mcB)^{<\infty}$ for all $m \geq 0$.
\end{defn}

\begin{lem}\label{L:highest_weight_cycles_are_cycles_infinite_strata}
The subcategory $\Perf(\mcX)^{\mathbf{O}\rm{-adm}} \subset \Perf(\mcX)$ is a thick stable subcategory that is closed under tensor products and symmetric powers, and $\rho_\ast((-) \otimes \mathbf{O}) : \QC(\mcX) \to \QC(\mcB)$ maps $\Perf(\mcX)^{\mathbf{O}\rm{-adm}}$ to $\DCoh(\mcB)$.
\end{lem}
\begin{proof}
	See \cite[Lem.~2.3]{halpernleistner2016equivariantverlindeformulamoduli} for the proof of the first claims. The definition of $\Perf(\mcX)^{\mathbf{O}\rm{-adm}}$ and \Cref{lem:perf_completion_infinite_strata} imply that $\rho_\ast((-)\otimes \mathbf{O})$ maps $\Perf(\mcX)^{\mathbf{O}\rm{-adm}}$ to $\DCoh(\mcB)$.
\end{proof}

\begin{prop}\label{P:nonabelian_localization_infinite_strata}
In the context of \Cref{D:perf_completion}, if $\mathbf{O} \in \Cycles(\mcX/\mcB)^{<\infty}$ and $P \in \Perf(\mcX)^{\mathbf{O}\rm{-adm}}$, then
\[
\chi(\mcX,P \otimes \mathbf{O}) = \sum_{\alpha} \chi\left(\mcZ_\alpha, \frac{[\tot_\alpha^\ast(P)]}{e(\bbN_{\mcZ_\alpha/\mcX})} \cdot \tot_\alpha^\sharp(\mathbf{O})\right),
\]
where the sum is well-defined because all but finitely many terms vanish for weight reasons.
\end{prop}
\begin{proof}
	As in the proof of \Cref{lem:perf_completion_infinite_strata}, restriction is an isomorphism $\rho_\ast(P\otimes \mathbf{O}) \cong (\rho|_{\mcX_{\leq \alpha}})_\ast(P \otimes \mathbf{O}|_{\mcX_{\leq \alpha}})$ for some $\alpha$. The stratification of $\mcX_{\leq \alpha}$ is finite because $\mcX_{\leq \alpha}$ is quasi-compact, so this follows from \eqref{E:perfect_formula}.
\end{proof}

\begin{ex}
    In the context of \Cref{L:perfect_highest_weight_cycles_quasi_smooth}, if $\mathbf{O} \in \Cycles(\mcX/\mcB)^{<\infty}$ is a homological shift of an invertible sheaf, then $\Perf(\mcX)^{\mathbf{O}\rm{-adm}}$ consists of $E \in \Perf(\mcX)$ such that for any fixed $m \geq 0$, one has
    \begin{equation}\label{E:weight_check_for_O_admissible}
        m \cdot \hwt(\ev_{1,\alpha}^\ast(E)) < \eta_\alpha - \wt(\mathbf{O}|_{\mcZ_\alpha})
    \end{equation}
    for all but finitely many $\alpha$, where $\eta_\alpha$ is given by \eqref{E:define_eta}.
\end{ex}

\subsection{Functoriality of highest weight cycles}\label{S:main_relative}

In this section we consider a commutative diagram
\begin{equation}\label{E:functoriality}
    \xymatrix{\mcX \ar[rr]^f \ar[dr]_{\rho_\mcX} & & \mcY \ar[dl]^{\rho_\mcY} \\ & \mcB & }
\end{equation}
and identify conditions under which $f_\ast : \QC(\mcX) \to \QC(\mcY)$ preserves categories of highest weight cycles.

\begin{ex}
    Pushforward along the morphism $\bbA_k^n / \bbG_m \to B\bbG_m$ corresponds to the functor taking a graded $k[x_1,\ldots,x_n]$-module (with each $x_i$ in degree $-1$) to the underlying graded vector space. $\Perf(\bbA^n_k/\bbG_m)_\beta^\wedge \cong \Cycles((\bbA^n/\bbG_m)/\Spec(k))^{<\infty}$ is the category of complexes of graded modules whose homology is finite dimensional in every weight space and vanishes for sufficiently large weight. It follows that pushforward maps $\Cycles((\bbA^n_k/\bbG_m)/\Spec(k))^{<\infty} \to \Cycles(B\bbG_m)^{<\infty}$
\end{ex}

Pre-composition with $(-)^n:\Theta \to \Theta$ for $n>0$ defines a morphism $\Filt_\mcR(\mcX) \to \Filt_\mcR(\mcY)$ that is an open and closed immersion \cite[Prop.~1.3.11]{Halpern-Leistner2014}, hence an isomorphism between connected components, and  it commutes with $\ev_1$. If you have a $\Theta$-stratification $\bigsqcup \mcS_i \subset \Filt_\mcR(\mcX)$, you can replace each $\mcS_i$ with its image under the morphism $(-)^{n_i}$ for some $n_i>0$ and obtain another $\Theta$-stratification without changing the image of the strata in $\mcX$. We call this operation ``scaling'' the stratification.

\begin{thm}\label{T:main_relative}
    Suppose that $\rho_\mcX$ and $\rho_\mcY$ satisfy the conditions of \Cref{D:perf_completion} and $f$ is quasi-compact and universally of finite cohomological dimension, and suppose that $\mcX = \bigcup_i \mcS_i$ and $\mcY = \bigcup_j \mcS'_j$ are (possibly infinite) $\Theta$-stratifications such that after scaling both stratifications, $\forall i$, $\exists j$ such that $f$ restricts to a $\Theta$-equivariant morphism $f : \mcS_i \to \mcS'_j$. Then $f_\ast : \QC(\mcX) \to \QC(\mcY)$ maps $\Cycles(\mcX/\mcB)^{<\infty}$ to $\Cycles(\mcY/\mcB)^{<\infty}$.
\end{thm}

The $\Theta$-equivariance condition means that either the $\Theta$-action on $\mcS'_j$ is trivial, or after regarding both strata as open substacks of $\Filt$, the morphism $\Filt_\mcR(\mcX) \to \Filt_\mcR(\mcY)$ maps $\mcS_i$ to $\mcS_j'$.

\begin{proof}
    We first prove the claim when both $\mcY$ and $\mcX$ are quasi-compact, hence the stratifications are finite. Using \Cref{L:local_cohomology_sod}, it suffices to prove the claim for cycles that are set theoretically supported on a single $\mcS_i$, and we can replace $\mcX$ with $\mcX_{\leq i}$ to assume that $\mcS_i$ is closed in $\mcX$. Let $\mcS'_j \subset \mcY$ be the unique stratum such that $f : \mcS_i \to \mcY$ factors through $\mcS'_j$. Using \Cref{P:functoriality_highest_weight_cycles} we may replace $\mcY$ with $\mcY_{\leq j}$ so that $\mcS'_j$ is closed.

    Now we claim if $F \in \Cycles(\mcX/\mcB\text{ on } \mcS_i)^{<v}$ then $f_\ast(F) \in \Cycles(\mcY/\mcB \text{ on } \mcS_i)^{<v}$, where the strata included in the notation indicates the subcategory of objects with the corresponding set theoretic support. We know that $F \cong \colim_{w \to -\infty} \radj{w}(F)$, $f_\ast(-)$ preserves filtered colimits, and $\QC_{\mcS'_j}(\mcY)^{<v}$ is closed under filtered colimits. Also, the $u^{th}$ associated graded object of this filtration of $F$ is of the form $(\ev_1)_\ast(G_u)$ for some $G_u \in \Cycles(\mcS_i/\mcB)^u$ by \Cref{L:constant_weight_cycles}. The morphism $f : \mcS_i \to \mcS'_j$ being $\Theta$-equivariant implies that the pullback $f^\ast$ preserves $\QC(-)^{\geq u}$ and $\QC(-)^{<u}$ \cite[Prop.~1.1.2(4)]{Halpern-Leistner2021DerivedConjecture}, so the right adjoint $f_\ast$ must preserve $\QC(-)^{<u}$. We know that $f_\ast$ preserves cycles automatically, so $f_\ast : \QC(\mcS_i) \to \QC(\mcS'_j)$ maps $\Cycles(\mcS_i/\mcB)^{<u}$ to $\Cycles(\mcS'_j/\mcB)^{<u}$. This implies that $f_\ast((\ev_1)_\ast(G_u)) \cong (\ev_1)_\ast ((f|_{\mcS_i})_\ast(G_u)) \in \Cycles(\mcY/\mcB \text{ on }\mcS'_j)^{<u}$ for all $u$. Combining all of these observations shows that $f_\ast(\radj{w}(F)) \in \Cycles(\mcY/\mcB \text{ on }\mcS'_j)^{<v}$ and $f_\ast(\ladj{w}(F)) \in \QC_{\mcS'_j}(\mcY)^{<w}$ for any $w \in \bbZ$, hence $\radj{w}(f_\ast(F)) \in \Cycles(\mcY/\mcB \text{ on }\mcS'_j)^{<v}$. \Cref{L:baric_truncation_preserves_highest_weight_cycles} then implies that $f_\ast(F) \in \Cycles(\mcY/\mcB \text{ on }\mcS'_j)^{<v}$.

    Now if the stratifications of $\mcY$ and $\mcX$ are infinite, the conditions of the theorem statement imply that $f^{-1}(\mcY_{\leq \alpha})$ is an open union of strata in $\mcX$ for any $\alpha$. It follows that for any $F \in \Cycles(\mcX/\mcB)^{<\infty}$, $f_\ast(F)|_{\mcY_{\leq \alpha}} \in \Cycles(\mcY_{\leq \alpha}/\mcB)^{<\infty}$. The other condition in \Cref{D:perf_completion} for $f_\ast(F)$ to lie in $\Cycles(\mcY/\mcB)^{<\infty}$ is that $\hwt(R\Gamma_{\mcS'_\alpha}(f_\ast(F)))<0$ for all but finitely many $\alpha$. This follows from the same condition on $F$, combined with the sharper claim proved in the last paragraph, that $f_\ast$ maps $\Cycles(\mcX/\mcB \text{ on } \mcS)^{<0}$ to $\Cycles(\mcY/\mcB \text{ on } \mcS')^{<0}$.
\end{proof}

\begin{ex}
    To see why one would want to scale the stratifications, consider a proper Deligne-Mumford stack $X$ with an action of $\bbG_m$. One can define this to mean there is an algebraic stack $\mcX$ with a morphism $\mcX \to B\bbG_m$ and an isomorphism $X \cong \rm{pt} \times_{B\bbG_m} \mcX$. The coarse moduli space $M$ of $X$ inherits a $\bbG_m$-action, and suppose $M$ admits an ample equivariant line bundle. The Bia{\l}ynicki-Birula stratification of $M/\bbG_m$ lifts to a $\Theta$-stratification of $\mcX$, but one might need to scale the strata in order to do this. As a result, if $\mcS_i$ is a stratum of $\mcX$, the projection $\mcS_i \to B\bbG_m$ might not be $\Theta$-equivariant for the $\Theta$-action on $B\bbG_m$ corresponding to the tautological cocharacter of $\bbG_m$, but this can be arranged after further scaling all of the strata. \Cref{T:main_relative} says that $f_\ast$ maps $\Cycles(\mcX)^{<\infty}$ to $\Cycles(B\bbG_m)^{<\infty}$ even without scaling the stratifications.
\end{ex}

\section{Applications}

\subsection{Wall crossing formulas}\label{S:wall_crossing}

Suppose that $\mcX$ is a quasi-smooth algebraic derived stack whose underlying classical stack admits a norm on graded points $\lVert-\rVert$ \cite[Def.~4.1.12]{Halpern-Leistner2014} and admits a good moduli space $\mcX \to X$ \cite{Al13} with $X$ proper over a field of characteristic $0$. Then for any $\ell$ in the relative N\'{e}ron-Severi group $\NS(\mcX)/\NS(X)$, or more generally any rational linear function on the component lattice of $\mcX$, the numerical invariant $\ell(f)/ \lVert f \rVert$ defines a finite $\Theta$-stratification $\mcX = \mcX^{\ell-\rm{ss}} \cup \bigcup_\alpha \mcS^\ell_\alpha$, for some indexing set $\alpha$. (See \cite[\S5.6]{Halpern-Leistner2014} for a complete discussion.)

The semistable locus $\mcX^{\ell\rm{-ss}}(\ell)$ admits a good moduli space that is projective over $X$ \cite[Thm.~5.6.1(2)]{Halpern-Leistner2014}, and choosing two values $\ell_1, \ell_2$ leads to a wall-crossing diagram
\[\xymatrix{\mcX^{\ell_1\rm{-ss}} \ar@{(=}[r] \ar[d]^{\rm{gms}} & \mcX \ar[d]^{\rm{gms}} & \mcX^{\ell_2\rm{-ss}} \ar[d]^{\rm{gms}} \ar@{(=}[l]\\
X_1 \ar[r] & X & \ar[l] X_2 }.\]
Applying the formula \eqref{E:quasi-smooth_localization_formula} to the $\Theta$-stratification for both values of $\ell$ gives a $K$-theoretic virtual wall-crossing formula for the index of any $F \in \Perf(\mcX)$,
\begin{equation}\label{E:wall_crossing}
\chi(\mcX^{\ell_1-\rm{ss}}, F) - \chi(\mcX^{\ell_2-\rm{ss}}, F) = \sum_\beta \chi\left(\mcZ^{\ell_2}_\beta,\frac{\tot_\beta^\ast(F)}{e(N_{\mcZ_\beta^{\ell_2}} \mcX)}\right) - \sum_\alpha \chi\left(\mcZ^{\ell_1}_\alpha,\frac{\tot_\alpha^\ast(F)}{e(N_{\mcZ_\alpha^{\ell_1}}\mcX)}\right).
\end{equation}
Here $\tot_\alpha : \mcZ_\alpha^{\ell_1} \to \mcX$ and $\tot_\beta : \mcZ_\beta^{\ell_2} \to \mcX$ are the centers of the unstable strata for $\ell_1$ and $\ell_2$ respectively.

This is somewhat simpler than wall-crossing formulas for cohomological integrals developed in \cite{Mochizuki2009} for two reasons: 1) One can define the $K$-theoretic integral on any stack with a proper good moduli space directly, whereas in cohomology one adds data to the moduli problem in order to obtain a DM stack, then shows that the integrals one gets are independent of these additional choices; and 2) One obtains a wall crossing formula without the auxiliary construction of a master space.

\begin{ex}[Sheaves on a surface]
	Let $S$ be a smooth projective surface, and let $H \in \NS(S)$ be an ample divisor class. Then the stack $\mcM^H_v$ of Gieseker semistable coherent sheaves on $S$ of class $v \in K_0^{\rm{num}}(S)$ is quasi-smooth and admits a projective good moduli space.
\end{ex}

\begin{ex}[Objects in a CY2 category]
	Let $\mcC$ be a smooth and proper $dg$-category $\mcC$ over a field $k$ of characteristic $0$, and assume that $\mcC$ is CY2 in the sense that the Serre functor on $\mcC$ is isomorphic to $[2]$. If $\sigma$ is a locally constant stability condition on $\mcC$ satisfying a mass-Hom bound \cite[\S2]{halpernleistner2025spaceaugmentedstabilityconditions}, then there is an algebraic derived stack $\mcM^\sigma_{v}$ of semistable objects in the heart of $\sigma$ and with class $v \in K_0^{\rm{num}}(\mcC)$. $\mcM^\sigma_v$ is a $0$-shifted symplectic derived stack, and hence it is quasi-smooth. It has affine diagonal, and its underlying classical stack admits a proper good moduli space \cite[Thm.~2.31]{halpernleistner2025spaceaugmentedstabilityconditions}. For any stability condition $\sigma'$ sufficiently close to $\sigma$, $\mcM_v^{\sigma'} \subset \mcM_v^\sigma$ is the semistable locus for a $\Theta$-stratification. This is shown in \cite[Prop.~4.4.5]{Halpern-Leistner2021DerivedConjecture} for derived categories of twisted K3 surfaces, the proof works for any CY2 category.
\end{ex}

\subsection{The stack of \texorpdfstring{$1$}{1}-dimensional sheaves on a surface.}\label{S:1D_sheaves}

There are elegant Verlinde-style formulas for the index of Atiyah-Bott $K$-theory classes \cite{Teleman2009TheCurve} on the stack of vector bundles of fixed rank on a smooth projective curve $C$. These were extended to index formulas on the stack of Higgs bundles in \cite{halpernleistner2016equivariantverlindeformulamoduli, andersen2017verlindeformulahiggsbundles}, and the stack of maps $\Map(C,X/G)$ for a linear representation $X$ in \cite{halpernleistner2023structuremoduligaugedmaps}. In this subsection we establish the foundational fact needed to formulate the analogous question for the stack of sheaves on a surface: that the cohomology of the Atiyah-Bott $K$-theory classes is finite dimensional.

Let $S$ be a smooth projective surface over a field, and let $H$ be an ample invertible sheaf on $S$. For any derived scheme or stack $T$, we will let $\pi : S \times T \to T$ denote the projection.

We will study the derived stack $\Coh_{1}(S)$ of flat families of coherent sheaves on $S$ whose support in every fiber has dimension $\leq 1$. This stack maps a derived scheme $T$ to the $\infty$-groupoid of objects $E \in \Perf(T \times S)$ that are $T$-flat, meaning that the functor $\QC(T) \to \QC(T \times S)$ taking $F \mapsto \pi^\ast(F) \otimes E$ is $t$-exact, and such that for any morphism $T' \to T$ from a classical scheme $T'$, $E|_{S_{T'}} \in \QC(S_{T'})^\heartsuit$ is a $T'$-flat coherent sheaf whose support in every fiber has dimension $\leq 1$. This is a locally finitely presented and quasi-smooth algebraic derived stack with affine diagonal \cite{Toen2004FromStacks}. We will let $\Coh_1^{\rm{pur}}(S) \subset \Coh_1(S)$ denote the open substack of points classifying pure sheaves, and let $\Coh_1^{\rm{pur}}(S)_r \subset \Coh_1^{\rm{pur}}(S)$ denote the open and closed substack of sheaves $E$ with $c_1(E) \cdot c_1(H)=r$.

\subsubsection*{The stratification}

For pure $1$-dimensional sheaves, the general theory of Gieseker semistability is equivalent to the usual theory of slope stability. We refer the reader to \cite{Huybrechts2010TheEdition} for a complete discussion, and simply recall the relevant facts here. For any point $[E] \in \Coh_1(S)$, we define the degree $\deg(E) := \chi(S,E)$, rank $\rank(E) := c_1(E) \cdot c_1(H)$, and slope $\mu(E):=\deg(E) / \rank(E)$. These are locally constant functions on $\Coh_1(S)$. A point $[E] \in \Coh_1^{\rm{pur}}(S)$ is said to be semistable if for all subsheaves $F \subset E$ with pure $1$-dimensional quotient, $\mu(F) \leq \mu(E)$. The substack of semistable sheaves of fixed degree and value of $c_1 \in \NS(S)_\bbQ$ is open and admits a projective good moduli space.

Every sheaf has a unique Harder-Narasimhan filtration, which is a filtration $E_q \subsetneq \cdots \subsetneq E_0 = E$ such that each $G_i := E_i/E_{i+1}$ is semistable and $\mu(G_0)<\cdots<\mu(G_q)$. If we fix $r = \rank(E)$ and equip every Harder-Narasimhan filtration with weights $w_i = M \cdot \mu(G_i)$, where $M = \rm{lcm}(1,\ldots,r)$ is chosen canonically to make $w_i$ an integer, then this identifies an open substack $\mcS \subset \Filt(\Coh_1^{\rm{pur}}(S)_r)$, which is a $\Theta$-stratification. The centers of this stratification are the moduli stacks of graded sheaves $G_0 \oplus \cdots \oplus G_q$ with each $G_i$ semistable of slope $w_i/M$ and with fixed values of $c_1(G_0),\ldots,c_1(G_q) \in \NS(S)$. Note that by labeling the strata this way, by the slopes $\mu_i$ and classes $c_1(G_i)$, the ranks and degrees of each $G_i$ are also determined. Each of these centers has a projective good moduli space.

\subsubsection*{Atiyah-Bott complexes}

One constructs perfect complexes on $\Coh_1(S)$ via an analog of the Atiyah-Bott construction. Consider the universal sheaf $E_{\rm{univ}}$ on $\Coh_1(S) \times S$, and let $p_2 : \Coh_1(S) \times S \to S$ and $\pi : \Coh_1(S) \times S \to \Coh_1(S)$ be the projections. Given a complex $F \in \Perf(S)$ and a partition $\lambda$, we define
\[
\mathbf{E}_\lambda^\ast(F) := R\pi_\ast(\bbS_\lambda(E_{\rm{univ}}) \otimes p_2^\ast(F)),
\]
where $\bbS_\lambda(-)$ is the extension to perfect complexes of the Schur functor associated to $\lambda$ \cite{Akin1982SchurComplexes}. We will also consider two line bundles
\begin{align*}
    \mcL_{\det} &:= \det(R\pi_\ast(E_{\rm{univ}}))^\ast \\
    \mcL_{\rm{rk}} &:= \det(R\pi_\ast(E_{\rm{univ}} \otimes p_2^\ast(H))) \otimes \mcL_{\det}
\end{align*}
The line bundle $\mcL_{\det}$ has weight $-\chi(E) = -\deg(E)$ with respect to the tautological action of $\bbG_m$ on $E_{\rm{univ}}$ by scaling, and $\mcL_{\rm{rk}}$ has weight $\rank(E) = c_1(H) \cdot c_1(E)$ with respect to this scaling $\bbG_m$-action.

\subsubsection*{Finiteness of cohomology}

\begin{thm}\label{T:1D_sheaves}
	For any $a, r>0$, $\mcL_{\det}^a \in \Cycles(\Coh_1^{\rm{pur}}(S)_r)^{<\infty}$, and for any $F \in \Perf(S)$, partition $\lambda$, and $b \in \bbZ$, $\mcL_{\rm{rk}}^b \otimes \mathbf{E}_\lambda^\ast(F)$ is $\mcL_{\det}^a$-admissible. In particular,\[
    \dim \left( \bigoplus_i H^i\left(\Coh_1^{\rm{pur}}(S)_r, \mcL_{\rm{det}}^a \otimes \mcL_{\rm{rk}}^b \otimes\mathbf{E}_\lambda^\ast(F) \right) \right) < \infty,\]
    where $\dim$ denotes dimension over the ground field of $S$.
\end{thm}
\begin{proof}
    A graded point $\gamma : B(\bbG_m)_k \to \Coh_1^{\rm{pur}}(S)_r$ corresponds to a direct sum decomposition $G = G_0 \oplus \cdots \oplus G_p$ of $G \in \Coh(S_k)$ along with a choice of integers $w_0 < \cdots < w_p$. We compute
    \begin{gather*}
        \wt(\gamma^\ast(\mcL_{\det})) = -\sum w_i \deg(G_i) = -\sum w_i \chi(S_k,G_i)  \\
        \wt(\gamma^\ast(\mcL_{\rm{rk}})) = \sum w_i \rank(G_i) = \sum w_i c_1(G_i) \cdot c_1(H)\\
        \hwt(\gamma^\ast(\mathbf{E}_\lambda^\ast(F))) \leq |\lambda| \cdot \max\{w_i\} \\
        \eta  = \wt(\det(\bbL_{\Coh_1(S),G}^+)) = \sum_{i>j} (w_i-w_j) \chi(S,G_j \otimes G_i^\ast [1]) = \sum_{i>j} (w_i-w_j) c_1(G_i) \cdot c_1(G_j).
    \end{gather*}
The last calculation uses the fact that $\bbL_{\Coh_1(S),G} \cong R\Hom(G,G[1])^\ast$ followed by Hirzebruch-Riemann-Roch.

We associate a stratum corresponding to HN filtrations whose associated graded objects have slopes $\mu_0 < \cdots < \mu_p$ to the vector $\vec{w} := M \cdot (\mu_0,\ldots,\mu_0,\mu_1,\ldots,\mu_1,\ldots) \in \bbZ^r$ where $\mu_i$ is repeated $\rank(E_i/E_{i+1})$-many times. Because there are only finitely many effective classes $D \in \NS(S)$ with $D \cdot c_1(H) \leq r$, there are only finitely many strata corresponding to any given vector $\vec{w}$.

By \Cref{L:perfect_highest_weight_cycles_quasi_smooth}, checking that $\mcL_{\det}^a \in \Cycles(\Coh_1^{\rm{pur}}(S)_r)$ amounts to checking that $\wt(\gamma_\alpha^\ast(\mcL_{\det})) < \eta_\alpha$ for the associated graded point of HN filtrations with all but finitely many values of $\vec{w} \in \bbZ^r$. Using the formula above, the weight of $\mcL_{\det}$ along the center of a stratum is $-\lVert \vec{w} \rVert^2$, where $\lVert-\rVert$ is the standard norm on $\bbZ^r$. On the other hand, because there are only finitely many possible values of $c_1(G_i)$ with $c_1(G_i) \cdot c_1(H) \leq r$, and therefore finitely many possible values of $c_1(G_i) \cdot c_1(G_j)$, there is a uniform upper bound $|\eta_\alpha| \leq K \lVert \vec{w} \rVert$, where $K$ depends on $r$ and the intersection pairing on the effective cone of $\NS(S)$. This implies that $\wt(\gamma^\ast(\mcL_{\det})) < \eta_\alpha$ for all but finitely many values of $\vec{w}$.

The fact that $\mcL_{\rm{rk}}^b \otimes \mathbf{E}_\lambda^\ast(F)$ is $\mcL_{\det}^a$-admissible is proved by verifying the inequality \eqref{E:weight_check_for_O_admissible} by an analogous calculation, observing that both $\wt(\gamma^\ast(\mcL_{\rm{rk}}))$ and $\hwt(\gamma^\ast(\mathbf{E}_\lambda^\ast(F)))$ have a uniform upper bound that is linear in $\lVert \vec{w} \rVert$ for the graded points of HN filtrations in $\Coh_1^{\rm{pur}}(S)_r$.

\end{proof}

\bibliography{localization}{}
\bibliographystyle{plain}

\end{document}